\documentclass{article}

\usepackage{amsmath, amsthm, amssymb}
\usepackage[all]{xy}

\newcommand{\CH}{C(H,\alpha)}
\newcommand{\HX}{H(X,\varphi)}
\newcommand{\HSig}{H(\Sigma,\sigma)}
\newcommand{\SX}{S(X,\varphi,P)}
\newcommand{\SSig}{S(\Sigma,\sigma,P)}
\newcommand{\UX}{U(X,\varphi,P)}
\newcommand{\USig}{U(\Sigma,\sigma,P)}
\newcommand{\GhX}{G^h(X,\varphi)}
\newcommand{\GsX}{G^s(X,\varphi,P)}
\newcommand{\GuX}{G^u(X,\varphi,P)}
\newcommand{\GhSig}{G^h(\Sigma,\sigma)}
\newcommand{\GsSig}{G^s(\Sigma,\sigma,P)}

\newtheorem{theo}{Theorem}[section]
\newtheorem{prop}[theo]{Proposition}
\newtheorem{cor}[theo]{Corollary}
\newtheorem{lemma}[theo]{Lemma}
\newtheorem{definition}[theo]{Definition}
\newtheorem{remark}[theo]{Remark}

\begin{document}


\title{Ring and module structures on dimension groups associated with a shift of finite type}
\author{D.B. KILLOUGH\footnote{Supported in part by an NSERC Scholarship}, \\
Department of Mathematics, Physics, and Engineering,\\
Mount Royal University,\\
Calgary, AB, Canada T3E 6K6,\\
bkillough@mtroyal.ca\\ 
I.F. PUTNAM \footnote{Supported in part by an NSERC Discovery Grant}, \\
Department of Mathematics and Statistics,\\
University of Victoria,\\
Victoria, B.C., Canada V8W 3R4,\\
putnam@math.uvic.ca\\
}

\maketitle

\begin{abstract}
We study invariants for shifts of finite type obtained as the K-theory of various $C^{\ast}$-algebras associated
with them. These invariants have been studied
intensely over the past thirty years since their introduction by  Wolfgang Krieger. They 
may be given quite concrete descriptions as inductive limits of simplicially ordered
free abelian groups. Shifts of finite type are special cases of Smale spaces and, in earlier work, the second author
has shown that the hyperbolic structure of the dynamics in a Smale space induces natural ring and module structures on
certain of these K-groups. Here, we restrict our attention to the special case of shifts of finite type and
obtain explicit descriptions in terms of the inductive limits.
\end{abstract}

\section{Introduction and Summary of Results} 

We provide a brief overview of our results. Precise versions of 
the definitions and results will follow in the later sections.

A Smale space, as defined by David Ruelle \cite{ruelle78}, is a compact metric space, $X$, together
with a homeomorphism, $\varphi$, which is hyperbolic. 
These include the basic sets of
Smale's Axiom A systems. Another special case of great interest are the shifts of finite type \cite{bowen78}, \cite{lindmarcus}
where the space, here usually denoted $\Sigma$, is the path space of a finite directed graph and 
the homeomorphism, $\sigma$, is the left shift.

Ruelle \cite{ruelle88} later showed how one may construct various $C^{\ast}$-algebras 
from a Smale space $(X, \varphi)$. These are best described as the groupoid $C^{\ast}$-algebras
associated with stable, unstable and homoclinic equivalence and we denote them 
by $\SX, \UX$ and $\HX$, respectively, where the first two depend on a choice of a set $P$ of periodic points.
It is worth mentioning that although the algebras depend on the set of periodic points, different 
choices will not alter the strong Morita equivalence class of the $C^{*}$-algebra, nor its
$K$-theory.
This extended earlier work of  Krieger \cite{krieger80} who considered
the case of shifts of finite type. For shifts of finite type, these $C^{\ast}$-algebras can be 
computed quite explicitly and are AF or approximately finite dimensional $C^{\ast}$-algebras.
As a consequence, their K-theory groups can be easily computed as inductive limits
of simplicially ordered free abelian groups. The invariant $K_{0}(\SSig)$ (or $K_{0}(\USig)$)
has been studied intensely and may be computed from the adjacency matrix of the underlying graph. 
This  was first  introduced by W. Krieger, who was motivated by Williams' notion of shift equivalence
for matrices. 
The group $K_{0}(\SSig)$ is usually called the \emph{dimension group} of $A$, as 
presented in section 7.5 of \cite{lindmarcus}. For a general Smale space, the computations of 
the $C^{*}$-algebras and their K-groups is a more complicated matter.

The original map, $\varphi$, in a Smale space provides automorphisms of the stable, 
unstable and homoclinic equivalence 
relations and hence automorphisms of the associated $C^{\ast}$-algebras, which we denote by
$\alpha$. In \cite{putnam96}, the second author showed that the hyperbolic nature 
of the dynamics  implies that the automorphism $\alpha$ of the homoclinic algebra $\HX$ is
asymptotically abelian. From this, various elements in the $E$-theory of Connes and Higson were constructed
in \cite{putnam96}. In doing so, it is necessary to pass from the discrete parameter ($\alpha$ is an action of the group 
of integers) to a continuous one. Thus, we are led to consider the mapping cylinder for 
the automorphism $(\HX, \alpha) $, which we denote by $\CH$, which now has an asymptotically abelian action
of the group $\mathbb{R}$. 
In consequence, the $K$-theory of  $\CH$ obtains a
natural product structure which makes it an ordered ring. For any $C^{*}$-algebra $A$, 
we write $K_{*}(A) = K_{0}(A) \oplus K_{1}(A)$, which is a $\mathbb{Z}_{2}$-graded abelian group.
(Since $\SSig$ and $\USig$ are AF-algebras, their $K_{1}$-groups are trivial.)
The ring structure is as a $\mathbb{Z}_{2}$-graded ring;
that is, the product
carries $K_{i}(\CH) \times K_{j}(\CH)$ to $K_{i+j}(\CH)$, where $i+j$ is interpreted mod 2.
 The $K$-groups of $\SX$ and $\UX$  are modules over this ring.
It is most convenient for us to regard  $K_{*}(\SX)$ as a right $K_{*}(\CH)$-module and 
$K_{*}(\UX)$ as a left $K_{*}(\CH)$-module.
It is worthwhile to note that the modules, are also all $\mathbb Z_{2}$-graded. 
Even in specific examples of more general 
Smale spaces where the K-groups may be computed, the computation of the  ring structure seems
quite difficult.

The main objective of this paper is to obtain explicit formulae for these ring and module structures for 
the special case of shifts of finite type. As a consequence, we find that the ring coincides with
the ring of continuous endomorphisms of the dimension group studied by Handelman \cite{handelman81, handelman09} in his analysis
of totally ordered dimension groups. In that sense, 
it is not new. But here it is being constructed from the $C^{*}$-algebra dynamical system $(\HX, \alpha) $
and is therefore closer in spirit to the program of noncommutative geometry. It seems possible
that further analysis of these $C^{*}$-algebras using finer tools of noncommutative geometry may
produce more dynamical invariants or information. It is also clear from these considerations that
any topological conjugacy between two shifts of finite type will induce
natural isomorphisms between the various $C^{*}$-algebras and their K-theories.

We present a summary of our results. Details will be provided in later sections. Let $A$ be a $K \times K$ matrix
with non-negative integer entries. Associated to $A$ is a finite directed graph with  vertex set 
$\{ 1, 2, \ldots, K \}$  and 
$A_{i,j}$ edges from vertex $i$ to vertex $j$, $ 1 \leq i, j \leq K$.
 There is also an associated shift of finite type, which we denote by
$(\Sigma, \sigma)$, the space $\Sigma$ consists of bi-infinite paths in the graph and $\sigma$ is the left shift.

The group $K_{0}( \SSig)$ is isomorphic to the inductive limit
$$
\xymatrix{
\mathbb{Z}^{K} \ar[r]^{v \mapsto vA} & \mathbb{Z}^{K} \ar[r]^{v \mapsto vA} & \mathbb{Z}^{K} \ar[r]^{v \mapsto vA} & \cdots .
}
$$
Here, we regard $\mathbb Z^{K}$ as being row vectors.
For a vector $v$ in $\mathbb{Z}^{K} $ and a non-negative integer $N$, we let $[v,N]$ denote the element in the inductive 
limit represented by $v$ as an element of the $N$th group in the sequence. The map $\sigma$ induces a 
natural automorphism of this group sending $[v,N]$ to $[vA, N]$. The group $K_{0}( \USig)$ may be obtained by replacing
$A$ by $A^{T}$, the transpose of $A$. Here, we prefer to instead consider $\mathbb Z^{K}$ as column vectors and
use multiplication by $A$ on the left instead of the right. 

The group $K_{0}(\HSig)$ is the inductive limit of the system 
$$
\xymatrix{
M_{K}({\mathbb Z}) \ \ar[r]^{X \mapsto AXA}  & \ M_{K}({\mathbb Z}) \ \ar[r]^{X \mapsto AXA} & \ M_{K}({\mathbb Z}) \ \ar[r] & \cdots
}
$$
and the automorphism induced by $\sigma$ sends $[X, N]$ to $[XA^{2}, N+1]$. In Krieger's original work, he used 
$\mathbb Z^{K} \otimes \mathbb Z^{K}$ instead of $M_{K}({\mathbb Z})$ and the map between these groups
was given by $A \otimes A^{T}$. We prefer the matrix notation as somewhat more convenient and better suited to
our computations of the product.

For a simple example, the three K-groups associated to the full $N$-shift (i.e. $A = [N ]$) are all 
$\mathbb Z[1/N]$. On the other hand, the automorphisms in the three cases are multiplication by $N$, by $\frac{1}{N}$, 
 and by $1$, respectively.

Finally, we consider the mapping cylinder $\CH$ and its K-groups (see section 2). We let $C(A)$ denote the centralizer of the 
matrix $A$ within the ring $M_{K}(\mathbb Z)$ and $B(A) = \{ AX - XA \mid X \in M_{K}(\mathbb Z) \}$. We prove that
$K_{0}(\CH)$ and $K_{1}(\CH)$ are the same inductive limits as above for $K_{0}(\HSig)$, but replacing 
$M_{K}(\mathbb Z)$ by the subgroup $C(A)$ and the quotient group $M_{K}(\mathbb Z)/B(A)$, 
respectively.
We provide  explanations of these
facts in section 3.

With this description, we can state our main results on the ring structure on $K_{*}(\CH)$.
\begin{theo}\label{K0XK0}
Let $(\Sigma, \sigma)$ be a mixing  SFT with $K \times K$ adjacency matrix $A$. 
Let $X, Y$ be in $C(A)$ and $N, M \geq 0$.  The product of
$[X,N]$ and $ [Y,M]$ in $K_0(\CH)$ is given by
$$
[X,N]\ast[Y,M] = [XY, N+M].
$$
\end{theo}

Without some additional structure such as we have indicated above, 
the K-theory of an arbitrary $C^{*}$-algebra does not possess a natural  ring structure.
The simplest case where this does occur is for $C(X)$, the $C^{*}$-algebra of continuous
functions on a compact Hausdorff space, $X$. Here, the key additional structure is that 
the $C^{*}$-algebra itself is commutative. In the end, the ring is
also commutative.
For our case, it 
is relatively easy to find an example of a matrix $A$ where the 
ring described above is non-commutative (and we provide one in section 3). To the best of our knowledge,
this is the of the $K$-theory of a $C^{*}$-algebra having a natural ring structure, where that ring
is not commutative.

The second case to consider is the product between elements of 
$K_0(\CH)$ and ones in $K_1(\CH)$.

\begin{theo}\label{K0XK1}
Let $(\Sigma, \sigma)$ be a mixing  SFT with $K \times K$ adjacency matrix $A$. 
Let $X$ be in $C(A)$, let $Y$ be in $M_{K}(\mathbb{Z})$ and let $M, N \geq 0$.
The product of $[X,N]$ in $K_0(\CH)$ and  $[Y + B(A), M]$ in $K_1(\CH)$ is
$$
[X,N] \ast [Y + B(A), M] = [XY + B(A), N + M ] \in K_1(\CH)
$$
and
$$
[Y + B(A), M] \ast [X,N] = [YX + B(A), N + M ] \in K_1(\CH).
$$
\end{theo}

Finally, we have the following case.
\begin{theo}\label{K1XK1}
Let $(\Sigma, \sigma)$ be a mixing  SFT with $K \times K$ adjacency matrix $A$. 
 For any $a, b$ in $K_{1}(\CH) $, their product  is zero in $K_{0}(\CH)$.
\end{theo}
We remark this result seems unlikely to persist in higher dimensional Smale spaces.

In \cite{putnam96}, the second author proved that integration against the measure of maximal entropy
 gave a trace on $\HX$ for any mixing Smale space $(X,\varphi)$, and that this easily extended to a 
trace on $\CH$.  Moreover, it was proved that the trace induced a ring homomorphism from $K_0(\CH)$ to $\mathbb{R}$. 
 In the case of a SFT we can write down explicitly what the trace is in terms of the inductive systems and
also the ring homomorphism.

\begin{cor}[Corollary to Thms. 3.3, 4.3 in \cite{putnam96}] \label{CHTrace}
Let $(\Sigma, \sigma)$ be a mixing  SFT  with $K \times K$ adjacency matrix $A$.
Let $\tau^{CH}$ be the trace on $\CH$.  Then for $[X,N] \in K_0(\CH)$, $\tau^{CH}_{\ast}[X,N] = \lambda^{-2N}u_lXu_r$. 
Where $\lambda$ is the Perron-Frobenius eigenvalue of $A$ and $u_l$, $u_r$ are the left and right Perron-Frobenius eigenvectors of $A$ normalized so that $u_lu_r = 1$.
\end{cor}

The next aim is to compute the structures of 
the usual dimension groups  as 
$K_{*}(\CH)$-modules. Since $\SSig$ and $\USig$ are both AF-algebras, their $K_{1}$-groups are trivial
and the only case we need to consider is the product between 
$K_{0}(\CH)$ and $K_{0}(\SSig)$ (and $K_{0}(\USig)$).

\begin{theo}\label{module}
Let $(\Sigma, \sigma)$ be a mixing  SFT with $K \times K$ adjacency matrix $A$. 
\begin{enumerate}
 \item 
For $v$ in $\mathbb{Z^K}$, considered as a row vector, $X$ in $C(A)$, $M, N \geq 0$, the product of $[v, N]$ in 
$K_{0}(\SSig)$ with $[X, M]$ in $K_{0}(\CH)$ is 
\[
  [v,N]\ast[X,M] = [vX,N+ 2M].
\]
 \item 
For $w$ in $\mathbb{Z^K}$, considered as a column vector, $X$ in $C(A)$, $M, N \geq 0$, the product of  $[X, M]$ in $K_{0}(\CH)$ with $[v, N]$ in 
$K_{0}(\USig)$ is 
\[
 [X,M]\ast[w,N] = [Xw,N + 2M].
\]
\end{enumerate}
\end{theo}

We observe the following.
\begin{cor}
The automorphisms $\alpha_{\ast}$ and $\alpha_{\ast}^{-1}$ on $K_0(\SSig)$ are induced by multiplication by $[A,0]$ and $[A,1]$ in $K_{0}(\CH)$,
 respectively.
\end{cor}

The existence of this module structure can be interpreted as a 
ring homomorphism from $K_{0}(\CH)$ into the endomorphism ring of the dimension group $K_{0}(\SSig)$.
The endomorphisms in the range of this homomorphism form a ring which was considered by Handelman \cite{handelman81, handelman82}, when $K_{0}(\SSig)$ has an order 
embedding with dense image in $\mathbb{R}^{k}$, for some $k$.
Moreover, these  endomorphisms are continuous in a sense described in \cite{handelman82}.
  A more general version also 
appears in \cite{handelman09}, where the integers have been replaced by a polynomial ring. 
In \cite{handelman82}, Handelman proves that this ring is part of  a complete invariant for the classification
of the group $K_{0}(\SSig)$ up to order isomorphism.

Here, our module structure gives a  map from $K_{0}(\CH)$ to the continuous endomorphisms  which is an isomorphism (this follows from the fact that the results of \cite{handelman09} and 
Theorems \ref{K0XK0} and \ref{module} yield exactly the same descriptions). The range of this isomorphism is precisely the ring of endomorphisms of the group $K_{0}(\SSig)$ (considered as a group without order) which commute with the automorphism $[v, N] \mapsto [vA, N]$ determined by the matrix $A$.

\begin{cor}\label{endo}
Let $(\Sigma, \sigma)$ be a mixing  SFT with $K \times K$ adjacency matrix $A$
and assume that there is an order embedding of $K_{0}(\SSig)$ in $\mathbb{R}^{k}$ 
(with the simplicial order) with dense image, for some $k$.  
The ring $K_{0}(\CH)$ is isomorphic to $End_{c}( K_{0}(\SSig))$, the ring of continuous endomorphisms of the dimension group
$ K_{0}(\SSig)$.
\end{cor}

The algebras $\SSig$ and $\USig$ have traces defined similarly to the trace on $\HSig$  in \cite{putnam96}. 
 We denote these by $\tau^s$ and $\tau^u$ respectively.  The maps on $K$-theory induced by these traces can be 
computed in terms of the inductive systems.

\begin{theo}
For $[v,N] \in K_0(\SSig)$, $[w,M] \in K_0(\USig)$ we have $\tau^s_{\ast}([v,N]) = \lambda^{-N}vu_r$, $\tau^u_{\ast}([w,M]) = \lambda^{-M}u_lw$. 
 Where $\lambda$ is the Perron-Frobenius eigenvalue of $A$, and $u_l$, $u_r$ are the left, right Perron-Frobenius eigenvectors of $A$.  
Moreover, if $[X,L] \in K_0(\CH)$, then we have
\begin{eqnarray*}
\tau^{CH}_{\ast}([v,N]\ast[X,L]) &=& \tau^s_{\ast}([v,N])\tau^{CH}_{\ast}([X,L]),\\ 
\tau^{CH}_{\ast}([X,L]\ast[w,M]) &=& \tau^{CH}_{\ast}([X,L])\tau^u_{\ast}([w,M]).
\end{eqnarray*}
\end{theo}

We next recall the definition of shift equivalence. 
The non-negative $n \times n$ integer matrix $A$ and the non-negative $m \times m$ integer matrix $B$ are said to be {\bf shift equivalent} if there exist non-negative integer matrices $R$ ($n \times m$) and $S$ ($m \times n$) and a positive integer $k$ such that:
$$
\begin{array}{ll}
\bullet \ RS = A^k  \hspace{25mm} &  \bullet \ AR = RB \hspace{25mm}\\
\bullet \ SR = B^k  \hspace{25mm} &  \bullet \ SA = BS. \hspace{25mm}
\end{array}
$$

Let $A$ and $B$ be two non-negative integral  matrices. We denote their associated shifts by $(\Sigma_{A}, \sigma_{A})$ and $(\Sigma_{B}, \sigma_{B})$, respectively.
We also choose sets $P_{A} \subset \Sigma_{A}$ and $P_{B}\subset \Sigma_{B}$ of periodic points. 
We say that the dimension triples $(K_0(\Sigma_{A}, \sigma_{A}, P_{A})),(K_0(\Sigma_{A}, \sigma_{A}, P_{A}))^+,(\alpha_{A})_{\ast})$
and $(K_0(\Sigma_{B}, \sigma_{B}, P_{B})),(K_0(\Sigma_{B}, \sigma_{B}, P_{B}))^+,(\alpha_{B})_{\ast})$ are isomorphic if there is 
an isomorphism between the groups, preserving the positive cones and intertwining the actions. This is equivalent to
the matrices $A$ and $B$  being shift equivalent as above. (See Theorem 7.5.8 of \cite{lindmarcus}.)

We claim that two SFTs with shift equivalent adjacency 
matrices also have isomorphic ring/module structures of $K$-theory.
This follows immediately from Corollary \ref{endo} and the last remark, but it 
is also fairly easy to give a direct proof 
following the ideas of Theorem 7.5.8 of \cite{lindmarcus} and we omit the details.

\begin{theo}
Let $(\Sigma_A, \sigma_A)$ and $(\Sigma_B, \sigma_B)$ be SFTs with adjacency matrices $A$, and $B$ 
respectively.  If $A$ and $B$ are shift equivalent, then
$$
\left(K_{\ast}(C(H_A,\alpha_A)), K_0(S_A), K_0(U_A)\right) \cong 
\left(K_{\ast}(C(H_B,\alpha_B)), K_0(S_B), K_0(U_B)\right)
$$
in the following sense: there exist
\begin{enumerate}
\item $\phi_H: K_{\ast}(C(H_A,\alpha_A)) \rightarrow K_{\ast}(C(H_B,\alpha_B))$  an isomorphism of ordered rings,
\item $\phi_S: K_0(S(\Sigma_A, \sigma_A)) \rightarrow K_0(S(\Sigma_B, \sigma_B))$ an isomorphism of ordered groups,
\item $\phi_U: K_0(U(\Sigma_A, \sigma_A)) \rightarrow K_0(U(\Sigma_B, \sigma_B))$ an isomorphism of ordered groups, and
\item for all $h \in K_{\ast}(C(H_A,\alpha_A))$, $s \in K_0(S(\Sigma_A, \sigma_A))$, $u \in K_0(U(\Sigma_A, \sigma_A))$ we have $\phi_S(s\ast h) = \phi_S(s)\ast\phi_H(h)$ and $\phi_U(h \ast u) = \phi_H(h)\ast\phi_U(u)$.
\end{enumerate}
\end{theo}

Our final result concerns the sense in which the two dimension groups
$K_{0}(\SSig)$ and $K_{0}(\USig)$ are dual. More concretely, 
we  look for a  result  of the form 
$Hom_R(K_0(\SSig),R) \cong K_0(\USig)$ for some subring $R$ of $K_0(\CH)$. 
As $K_0(\CH)$ is, in general, non-commutative, the subring $R$ should lie in the center of the ring $K_0(\CH)$.

 Let $R_{A}$ be the subring of $K_{0}(\CH)$ generated by $[A, 0]$ and $[A,1]$.
It follows easily from Theorem \ref{K0XK0} that these elements are inverses of each other and, 
in particular, $R_{A}$ contains the unit of $K_{0}(\CH)$.
Moreover, 
 it is clear that $R_{A}$ is contained in the center of the ring $K_0(\CH)$, $Z(K_0(\CH))$.  
 In many, but not all, cases $R_{A} = Z(K_0(\CH))$. We present an example in section \ref{secDuality}
 where they differ. We note the following standard  description of $R_{A}$.

\begin{prop}
 Let the minimal polynomial of $A$ be $m_A(x) = x^l(x^k + a_{k-1}x^{k-1} + \cdots + a_0)$ with $a_0 \neq 0$. 
 Then
$$
R_{A} \cong \mathbb{Z}[x,x^{-1}]/<p_A(x)>,
$$
where $p_A(x) = x^k + a_{k-1}x^{k-1} + \cdots + a_0$.
\end{prop}
 
 Our duality result is the following.

\begin{theo}\label{ModDuality}
Let $(\Sigma, \sigma)$ be a mixing SFT, and $R_{A}$ the subring of $K_0(\CH)$ generated by 
$[A, 0]$ and $[A,1]$. Then
$$
Hom_{R_{A}}(K_0(\SSig),R_{A}) \cong K_0(\USig)
$$ 
as left $R_{A}$-modules.
\end{theo}

In section 2, we provide background information on Smale spaces, their $C^{*}$-algebras and their $K$-theory. We turn our attention to the 
special case of shifts of finite type in section 3 and we give proofs of the formulae for product from Theorems \ref{K0XK0}, \ref{K0XK1}
and \ref{K1XK1}.
Section 4 deals with the duality result.

The alert reader will have noticed that all of our results above assume that the shift of finite type is mixing. In section 5,
we consider the irreducible case and give descriptions of the $C^{*}$-algebras and their K-theory which effectively allows
an extension of our main results to this case.

It is a pleasure to thank the referees for many helpful comments.

\section{Smale spaces}

\subsection{Dynamics}

Smale spaces were defined by Ruelle in \cite{ruelle78}, based on the Axiom A systems studied by Smale in \cite{smale67}. For a precise definition and many results relevant to this paper, we refer the reader to \cite{putnam96}. 

Roughly speaking, a Smale space is a topological dynamical system $(X,\varphi)$ in which $X$ is a compact metric space with distance function $d$, and $\varphi$ is a homeomorphism.  The structure of $(X,\varphi)$ is such that each point $x \in X$ has two local sets associated to it: A set, $V^s(x,\epsilon)$, on which the map $\varphi$ is (exponentially) contracting; and a set, $V^u(x,\epsilon)$, on which the map $\varphi^{-1}$ is contracting.  We call these sets the local stable and unstable sets for $x$. Furthermore, $x$ has a neighbourhood, $U_x$ that is isomorphic to $V^u(x,\epsilon) \times V^s(x,\epsilon)$. In other words, the sets $V^u(x,\epsilon)$ and $V^s(x,\epsilon)$ provide a coordinate system for $U_x$ such that, under application of the map $\varphi$, one coordinate contracts, and the other expands.  We denote this homeomorphism by $[\cdot,\cdot] : V^u(x,\epsilon) \times V^s(x,\epsilon) \rightarrow U_x$. 

We now define three equivalence relations on $X$.
  We say $x$ and $y$ are \emph{stably equivalent} and write $x \stackrel{s}{\sim} y$ if 
$$
\lim_{n \rightarrow +\infty}d(\varphi^n(x),\varphi^n(y)) = 0.
$$
We say $x$ and $y$ are \emph{unstably equivalent} and write $x \stackrel{u}{\sim} y$ if 
$$
\lim_{n \rightarrow -\infty}d(\varphi^n(x),\varphi^n(y)) = 0.
$$
Finally, we say $x$ and $y$ are \emph{homoclinic} and write $x \stackrel{h}{\sim} y$ if $x \stackrel{s}{\sim} y$ and $x \stackrel{u}{\sim} y$.  We denote these three equivalence relations by $G^s$,$G^u$, and $G^h$.
We also let $V^{s}(x), V^{u}(x)$ and $V^{h}(x)$ denote the 
equivalence classes of a point $x$ in each.

Throughout the first sections of this paper, we will assume our Smale space is mixing: for every pair of
non-empty open sets, $U$ and $V$, there is a constant $N$ such that $\varphi^{n}(U) \cap V$ is non-empty if
$n \geq N$.

\subsection{$C^{\ast}$-algebras}

In general, an equivalence relation is a  groupoid and we construct groupoid $C^{\ast}$-algebras from  stable, unstable, and homoclinic equivalence.  The construction of these $C^{\ast}$-algebras for a given Smale space is  originally due to Ruelle (\cite{ruelle88}).  In the case of a shift of finite type, these are the algebras studied by Cuntz and Krieger in \cite{cuntzkrieger}, \cite{krieger80}.   We summarize as follows.

The groupoid for homoclinic equivalence may be endowed with an \'{e}tale topology, as in \cite{putnam96}. In \cite{putnamspielberg},
it was shown to be amenable. We let 
$\HX$ denote its $C^{\ast}$-algebra.

For stable and unstable equivalence, we proceed as  in \cite{putnamspielberg}.  We first fix a finite $\varphi$-invariant subset of $X$, $P$. We note that in an irreducible Smale space, the set of periodic points is dense, so there are plenty of choices. 
We then consider the set of all points in $X$ that are unstably equivalent to a point in $P$, call this $V^u(P)$.  This set may be endowed with a natural topology in which it
is locally compact and Hausdorff. The groupoid that we actually use to construct our stable algebra is then the groupoid of stable equivalence restricted to the set $V^u(P)$. As described in \cite{putnamspielberg}, this groupoid also has a natural \'{e}tale topology and is amenable. We let $\SX$ denote the reduced $C^{\ast}$-algebra. For different choices
of $P$, the groupoids are equivalent to each other in the sense of Muhly, Renault and Williams 
\cite{muhlyrenaultwilliams}, and the groupoid $C^{\ast}$-algebras are Morita equivalent (see \cite{putnamspielberg}).  The construction of $\UX$ is similar.

In \cite{putnam96}, it is shown that $\HX$ is contained in the multiplier algebras of both 
$\SX$ and $\UX$; that is there are products $\HX \times \SX, \SX \times \HX \rightarrow \SX$ and similarly for
$\UX$. Although the definition of $\SX$ and $\UX$ are slightly different in \cite{putnam96}, the same
formulae still yield well-defined products.

The homeomorphism $\varphi$ yields a $\ast$-automorphism on each of the three algebras associated to $(X,\varphi)$. For $f$,
a continuous function of compact support on any one of our groupoids, 
we define $\alpha(f)$ by $\alpha(f)(x,y) = f(\varphi^{-1}(x),\varphi^{-1}(y))$.  There are several asymptotic commutation results that arise from $\alpha$.  We recall the following theorem from \cite{putnam96}.
Let and $a,b$ be in $\HX$, let $c$ be in $\SX$ and  $d $ be in $\UX$. Then we have
\begin{enumerate}
 \item $||[\alpha^n(a), b]|| \rightarrow 0$ as $n \rightarrow \pm \infty$, 
\item $||[\alpha^n(a), c]|| \rightarrow 0$ as $n \rightarrow - \infty$, 
\item $||[\alpha^n(a), d]|| \rightarrow 0$ as $n \rightarrow + \infty$.
\end{enumerate}

Let 
$$\CH = \{ f: \mathbb{R} \rightarrow \HX \ | \ f \ \textrm{cont.}, \ f(t+1) = \alpha(f(t)) \}$$
 and 
define  automorphisms $\alpha_t(f)(s) = f(t+s), s \in \mathbb{R}$, 
for each $t$ in $\mathbb{R}$.  First, we note that, for $f$ in $\CH$ and $c$ in $\SX$, 
the formulae $f \times c \rightarrow f(0)c$ and $c \times f \rightarrow c f(0)$ also 
define products from $\CH \times \SX, \SX \times CH \rightarrow \SX$ and similarly for $\UX$.
Secondly, it is an immediate consequence of the result above that, for
$f, g$  in $\CH$,  $c$  in $\SX$ and  $d $  in $\UX$, we have
\begin{enumerate}
 \item $||[\alpha_{t}(f), g]|| \rightarrow 0$ as $t \rightarrow \pm \infty$, 
\item $||[\alpha_{t}(f), c]|| \rightarrow 0$ as $t \rightarrow - \infty$, 
\item $||[\alpha_{t}(f), d]|| \rightarrow 0$ as $t \rightarrow + \infty$.
\end{enumerate}

\subsection{K-Theory}\label{KCHsection}

In passing from $\HX$ to $\CH$, we  have turned the integer parameter to a real one. The significance 
is that these maps now define asymptotic morphisms in the sense of Connes and Higson, \cite{conneshigson}.
Specifically, the map $\CH \otimes \CH \rightarrow \CH$ given by $f \otimes g \mapsto \alpha_t(f)\alpha_{-t}(g)$ 
is an asymptotic morphism which determines an element in the group $E(\CH \otimes \CH, \CH)$.  
This in turn yields a map $K_{\ast}(\CH \otimes \CH) \rightarrow K_{\ast}(\CH)$.
The idea is simple: if $\rho_{t}$ is an asymptotic homomorphism from a $C^{\ast}$-algebra $A$ to another $B$, and 
$p$ is a projection in $A$, then $\rho_{t}(p)$ is almost a projection in $B$, for large values of $t$. More precisely, 
$\rho_{t}(p)^{2} - \rho_{t}(p)$ and $\rho_{t}(p)^{\ast} - \rho_{t}(p)$ are small in norm and so the spectrum of 
$(\rho_{t}(p) + \rho_{t}(p)^{\ast})/2$ is concentrated near $0$ and $1$. We may use functional calculus with the 
function $\chi_{(1/2,\infty)}$ to obtain a projection in $B$. This depends on $t$, but the projections vary
continuously and hence determine a unique element in the $K$-theory of $B$.

Combining this with a Kunneth Theorem (e.g., Theorem 23.1.3 in \cite{blackadar}) gives a map 
$K_{\ast}(\CH) \otimes K_{\ast}(\CH) \rightarrow K_{\ast}(\CH)$. That is, a 
ring structure on the group $K_{\ast}(\CH)$.
Similarly, the asymptotic morphisms  $\SX \otimes \CH \! \rightarrow \! \SX$ defined by
 $c \otimes f \rightarrow c f(-t)$ and
$\UX \otimes \CH \! \rightarrow \! \UX$ defined by $f \otimes d \rightarrow f(t)d$
give rise to right and left $\CH$-module 
structures for $\SX$ and $\UX$, respectively.

We now define the ring structure on $K_0(\CH)$ in more detail. For any
$C^{\ast}$-algebra $A$ and $n \geq 1$, we let $M_{n}(A)$ denote the $C^{\ast}$-algebra of $n \times n$ matrices over $A$.
We also let $P_{n}(A)$ denote the set of projections in $M_{n}(A)$. We let $P_{\infty}(A)$ denote the union
of all $P_{n}(A)$.
For $f \in M_n(\CH)$, $g \in M_m(\CH)$ we define $(f \times g)_t \in M_{nm}(\CH)$ componentwise by
$$
((f \times g)_t)_{(i,j)(i',j')} = (f_{i,i'} \times g_{j,j'})_t = \alpha_t(f_{i,i'})\alpha_{-t}(g_{j,j'}),
$$
for $1 \leq i, i' \leq n, 1 \leq j, j' \leq m, t \geq 0$.
For $p \in P_n(\CH)$, $q \in P_m(\CH)$, there exists $T >0$ such that

$$
\chi_{(1/2,\infty)}(p \times q)_t \in P_{nm}(\CH) 
$$
for $t \geq T$.
The function $t \mapsto \chi_{(1/2,\infty)}(p \times q)_t$ is continuous, so for $t \geq T$,  $\chi_{(1/2,\infty)}(p \times q)_t$ 
forms a continuous path of projections in $P_{nm}(CH)$.  Thus, for $t_1, t_2 \geq T$
$$
\left[\chi_{(1/2,\infty)}(p \times q)_{t_1}\right]_0 = \left[\chi_{(1/2,\infty)}(p \times q)_{t_2}\right]_0.
$$
The following proposition gives a concrete form for the product on $K_0(\CH)$.  We state the result without proof.

\begin{prop}\label{Kprod}
For $p,q \in P_{\infty}(\CH)$ we define the product $[p]_0[q]_0 \in K_0(\CH)$ by 
$$
[p]_0[q]_0 = \lim_{t \rightarrow \infty} \left[\chi_{(1/2,\infty)}(p \times q)_t \right]_0
$$
\end{prop}

Our product will actually be defined on $K_{\ast}(\CH) = K_{0}(\CH) \oplus K_{1}(\CH)$, as a $\mathbb{Z}_2$-graded group. To do so, we simply use the fact 
that $K_{1}(\CH) \cong K_{0}(S \CH)$, where $SA$ denotes the suspension of a $C^{\ast}$-algebra $A$. Then an obvious extension yields asymptotic morphisms
 $\CH \otimes S \CH \rightarrow S \CH$, $S \CH \otimes  \CH \rightarrow S \CH$ and $S \CH \otimes  S \CH \rightarrow S^{2} \CH$.

As a first step toward computing the product, we describe the $K$-theory of the mapping cylinder, $\CH$.

We have the following short exact sequence.
$$ 0 \longrightarrow SH \stackrel{\iota}{\longrightarrow} \CH \stackrel{e_0}{\longrightarrow} \HX \longrightarrow 0 $$
The map $e_0$ is evaluation at $0$, and $\iota(f)(s) = \alpha^k(f(s-k))$ for $k \leq s \leq k+1$.
We thus get the following 6-term exact sequence of $K$ groups.
$$
\xymatrix{
K_0(SH) \ar[r]^-{\iota_{\ast}} & K_0(\CH)\ar[r]^-{(e_0)_{\ast}} & K_0(\HX)\ar[d] \\
K_1(\HX) \ar[u] & K_1(\CH) \ar[l]^-{(e_0)_{\ast}} & K_1(SH)\ar[l]^-{\iota_{\ast}}
}
$$
For a general irreducible Smale space, it is a difficult problem to obtain a nice description of $K_{\ast}(\HX)$, and hence difficult to describe 
concretely $K_{\ast}(\CH)$.  However, in section \ref{SFTsection} we will see that in the SFT case, the above 6-term exact sequence 
can be used to describe $K_{\ast}(\CH)$ in a concrete way.

\section{Shifts of finite type}\label{SFTsection}

\subsection{Dynamics}
We begin with a finite connected directed graph $G$, with no sources or sinks. That is, $G$ has a vertex set $V(G)$, an edge set $E(G)$ and initial 
and terminal maps $i, t: E(G) \rightarrow V(G)$.  We let $A$ denote the (non-degenerate) adjacency matrix of the graph. That is, if we enumerate the vertices of $V(G) = \{ v_{1}, v_{2}, \ldots, v_{K} \}$, then $A$ is the $K \times K$ matrix whose
$i,j$ entry is the number of edges $e$ with $i(e)=v_{i}$ and $t(e) = v_{j}$. 

We let $\Sigma$ be the associated bi-infinite path space of the graph:
\[
 \Sigma = \{ (e_{n})_{n \in \mathbb{Z}} \mid e_{n} \in E(G), t(e_{n}) = i(e_{n+1}), \text{ for all } n \in \mathbb{Z} \}.
\]
We also let $\sigma$ denote the left shift map. See \cite{lindmarcus} for a complete treatment. The system 
$(\Sigma, \sigma)$ is a Smale space \cite{putnam09} and is mixing precisely when the matrix $A$ is primitive; i.e. 
there is a positive integer $M$ such that $A^{N}$ has no zero entries for all $N \geq M$.

\subsection{$C^{\ast}$-Algebras}\label{C*SFT}

  We begin by describing the algebras $\HSig$, and $\SSig$.  These constructions follow (with slight notational modifications) those in \cite{putnam09}.  We present the details here for completeness, and because some of our later proofs will use this notation.  With $A$ primitive and $M$ as above,
fix $N \geq M$, $v_i, v_j \in V(G)$.  Define
$$ \Xi_{N,v_i,v_j} = \{\xi=(\xi_{-N+1}, \cdots, \xi_N) \ \ |\ \ t(\xi_N)=v_j, \ i(\xi_{-N+1})=v_i, t(\xi_{n}) = i(\xi_{n+1}), \forall -N < n < N  \}.$$
Notice that $\Xi_{N,v_i,v_j}$ consists of all allowable paths of length $2N$ in $G$ with initial vertex $v_i$ and terminal vertex $v_j$, so $\# \Xi_{N,v_i,v_j} = A^{2N}_{ij} > 0$.
For $\xi \in \Xi_{N,v_i,v_j}$ define
$$V_{N,v_i,v_j}(\xi)=\{x \in \Sigma \ \ |\ \  x_n = \xi_n \ \ \  \forall -N+1 \leq n \leq N \}.$$
Note that for fixed $N$, $V_{N,v_i,v_j}(\xi)$ and $V_{N,v_i',v_j'}(\eta)$ intersect only if $\xi = \eta$, $v_i = v_i'$, and $v_j = v_j'$.  Now let $\xi \ ,\eta \in \Xi_{N,v_i,v_j}$. Define
$$
E_{N,v_i,v_j}(\xi,\eta) \! = \! \{ (x,y) \ | \ \{x_n\}_{-N+1}^{N} = \xi, \ \{y_n\}_{-N+1}^{N} = \eta, \ x_n = y_n \ \forall n > N, \ n < -N+1 \}.
$$
Then
\begin{enumerate}
\item $E_{N,v_i,v_j}(\xi,\eta)$ and $E_{N,v_i',v_j'}(\xi',\eta')$ intersect only if $\xi = \xi'$, $\eta = \eta'$, $v_i = v_i'$, $v_j = v_j'$.
\item The sets $E_{N,v_i,v_j}(\xi,\eta)$ for $N \geq 1$, $v_i, \ v_j \in V(G)$, $\xi, \ \eta \in \Xi_{N,v_i,v_j}$ form a neighbourhood base of compact, open sets for the topology on $G^h(\Sigma, \sigma)$.
\end{enumerate}
Now let 
$e_{N,v_i,v_j}(\xi, \eta) = \chi_{E_{N,v_i,v_j}(\xi,\eta)} \in C_c(G^h(\Sigma, \sigma)).$  Note that 
$ \overline{span \{ e_{N,v_i,v_j}(\xi, \eta) \} } = \HSig$.
Consider the product of two such functions: 
$$e_{N,v_i,v_j}(\xi, \eta) \ast e_{N,v_i',v_j'}(\xi', \eta')(x,y) = \sum_{x \stackrel{h}{\sim z}}e_{N,v_i,v_j}(\xi, \eta)(x,z)e_{N,v_i',v_j'}(\xi', \eta')(z,y).$$
It is easily verified that this reduces to
\begin{displaymath}
e_{N,v_i,v_j}(\xi, \eta) \ast e_{N,v_i',v_j'}(\xi', \eta') = \left\{ \begin{array}{ll}
e_{N,v_i,v_j}(\xi, \eta') & \textrm{if $\eta = \xi'$} \\
0 & \textrm{otherwise}\\
\end{array} \right.
\end{displaymath}
Now let 
$$
H_{N,v_i,v_j} = span \{e_{N,v_i,v_j}(\xi, \eta) \ | \ \xi, \, \eta \in \Xi_{N,v_i,v_j} \}\cong M_{k(N,v_i,v_j)}({\mathbb C}).
$$
where $k(N,v_i,v_j) = \# \Xi_{N,v_i,v_j} = A^{2N}_{v_1,v_2}$.  
Now we define
$$
H_N = span(\{e_{N,v_i,v_j}(\xi, \eta) \ | \ \xi, \, \eta \in \Xi_{N,v_i,v_j}; \ v_i, \, v_j \in V(G)\}),
$$
and notice that
$$
H_N = \bigoplus_{v_i \in V(G)} \bigoplus_{v_j \in V(G)} H_{N,v_i,v_j} = \bigoplus_{(v_i,v_j) \in V(G) \times V(G)} H_{N,v_i,v_j} \cong \bigoplus_{(v_i,v_j)}M_{k(N,v_i,v_j)}({\mathbb C}).
$$

Now $H_N \subset H_{N+1}$, and $\HSig$ is the direct limit of the $H_N$'s. To see how $H_N$ is embedded in $H_{N+1}$,
observe that
$$
e_{N,v_i,v_j}(\xi, \eta) = \sum_{y_1 \in E_i}\sum_{y_2 \in E_j}e_{N+1,v_l,v_k}(y_1\xi y_2, y_1\eta y_2),
$$
where $i(y_1) = v_l$, $t(y_2) = v_k$, $E_i = \{ y \in E(G) \ | \ t(y) = v_i \}$, and $E_j = \{ y \in E(G) \ | \ i(y) = v_j \}$. In particular, $H_{N+1,v_l,v_k}$ contains $A_{li}A_{jk}$ copies of $H_{N,v_i,v_j}$.

We now describe the action of $\alpha$ on $\HSig$. 
$$
\alpha(e_{N,v_i,v_j}(\xi, \eta)) = \sum_{k}\sum_{\xi' \in \Xi_{1,v_j,v_k}}e_{N+1,v_i,v_k}(\xi \xi',\eta \xi'),
$$
and
$$
\alpha^{-1}(e_{N,v_i,v_j}(\xi, \eta)) = \sum_{l}\sum_{\xi' \in \Xi_{1,v_l,v_i}}e_{N+1,v_l,v_j}(\xi' \xi,\xi' \eta).
$$
In particular $\alpha$ and $\alpha^{-1}$ map $H_N$ into $H_{N+1}$.
  
The construction of $\SSig$ is very similar. We briefly outline the details. Fix a finite $\sigma$-invariant set $P \subset \Sigma$.
Fix $N \geq M$, $v_i \in V(G)$.  Define
$$ 
\Xi_{N,v_i} = \{\xi=(\xi_{-N+1}, \cdots, \xi_N) \ \ |\ \ t(\xi_N)=v_i, \ i(\xi_{-N+1}) = i(p_{-N}) \ \textrm{for some} \ p \in P \}.
$$
Again we mention that $\Xi_{N,v_i}$ is non-empty, as $A^{2N}$ is strictly positive.
For $\xi \in \Xi_{N,v_i}$ we can extend $\xi$ backwards by setting $\xi_{-n} = p_{-n}$ for $n>N-1$.
Now for $\xi \in \Xi_{N,v_i}$ we define
$$
V_{N,v_i}(\xi)=\{x \in \Sigma \ \ |\ \  x_n = \xi_n \ \ \  \forall n \leq N \}.
$$
Note that for fixed $N$, $V_{N,v_i}(\xi)$ and $V_{N,v_j}(\eta)$ intersect only if $\xi = \eta$, and $v_i = v_j$.  Now let $\xi \ ,\eta \in \Xi_{N,v_i}$. Define
$$
E_{N,v_i}(\xi,\eta) = \{ (x,y) \ | \ x_n = \xi_n, \ y_n = \eta_n \ \forall n \leq N, \ x_n = y_n \ \forall n > N \}.
$$
The collection of sets $\{E_{N,v_i}(\xi,\eta) \}$ forms a clopen base for the topology on $\GsSig$, and we are left to consider functions of the form
$$
e_{N,v_i}(\xi,\eta) = \chi_{E_{N,v_i}(\xi,\eta)}.
$$
Proceeding as we did above for $\HSig$, we see that for fixed $N$ and $i$
$$
e_{N,v_i}(\xi,\eta)\ast e_{N,v_i}(\xi',\eta') = \left \{ \begin{array}{cl}
e_{N,v_i}(\xi,\eta') & \textrm{if $\eta = \xi'$}, \\
0 & \textrm{if $\eta \neq \xi'$}.
\end{array} \right.
$$
As above, we let $S_{N,v_i} = span\{e_{N,v_i}(\xi,\eta) \ | \  \xi,\eta \in \Xi_{N,v_i} \}$ and notice that
$$
S_{N,v_i} \cong M_{k(N,v_i)}(\mathbb{C}),
$$
where $k(N,v_i)$ is the number of paths of length $2N$ starting at a vertex of $p \in P$ and ending at $v_i$. 
$$
S_N = \bigoplus_{v_i \in V(G)} S_{N,v_i} \cong \bigoplus M_{k(N,v_i)}(\mathbb{C}).
$$
Finally we notice that $S_N \subset S_{N+1}$ and let $\SSig$ be the direct limit of the $S_N$'s. Similar to the above,
$$
e_{N,v_i}(\xi,\eta) = \sum_{y \in S}e_{N+1,v_k}(\xi y,\eta y),
$$
where $t(y) = v_k$ and $S = \{ y \in E(G) \ | \ i(y) = v_i \}$. So we see that $S_{N,v_k}$ contains $A_{ik}$ copies of $S_{N,v_i}$.

Similar to the $\HSig$ case we see that
$$
\alpha(e_{N,v_i}(\xi,\eta)) = \sum_{k}\sum_{\xi' \in \Xi_{1,v_i,v_k}}e_{N+1,v_k}(\xi \xi',\eta \xi'),
$$
and
$$
\alpha^{-1}(e_{N,v_i,v_j}(\xi, \eta)) = e_{N+1,v_i,v_j}(\xi,\eta).
$$

\subsection{$K$-Theory}\label{K(CH)}
We now describe the $K$-theory for $\HSig$ and $\SSig$ in the case that $(\Sigma, \sigma)$ is mixing. The irreducible case is handled in section \ref{irred}. The computation of the $K$-theory follows easily from the description as AF-algebras in section \ref{C*SFT}.  These results also appear in \cite{putnam96}.  We begin with $\HSig$.  This was first computed by Krieger in \cite{krieger80}.

As $\HSig$ is AF, $K_1(\HSig) = 0$.  $K_0(H_N) = \bigoplus_{V\times V}\mathbb{Z} \cong M_{K}({\mathbb Z})$, and $K_0(\HSig)$ is the inductive limit of the following system.
$$
\xymatrix{
M_{K}({\mathbb Z}) \ \ar[r]^{X \mapsto AXA}  & \ M_{K}({\mathbb Z}) \ \ar[r]^{X \mapsto AXA} & \ M_{K}({\mathbb Z}) \ \ar[r] & \cdots
}
$$
The inductive limit group can be described as
$$K_0(H) \cong (M_{K}({\mathbb Z}) \times {\mathbb N})/\sim.$$
Where, for $n \leq k$, $(X,n) \sim (Y,k)$ if and only if $A^{k-n+l}XA^{k-n+l} = A^lYA^l$ for some $l \in {\mathbb N}$. We denote the equivalence class of $(X,N)$ under $\sim$ by $[X,N]$.

Recall the automorphism $\alpha: H \rightarrow H$.  We now wish to describe $\alpha_{\ast}:K_0(H) \rightarrow K_0(H)$.  Referring back to section \ref{C*SFT}, it is another straightforward calculation to see that for $[X,N] \in  K_0(H)$ we have $\alpha([X,N]) = [XA^2,N+1]$.  Similarly, $\alpha^{-1}([X,N]) = [A^2X,N+1]$.

We now briefly outline the computation of $K_{\ast}(\SSig)$.
As in the case of $\HSig$, $\SSig$ is AF and hence $K_1(\SSig) = 0$, and $K_0(\SSig)$ is the direct limit of the following system.
$$
\xymatrix{
\mathbb{Z}^{K} \ar[r]^{v \mapsto vA} & \mathbb{Z}^{K} \ar[r]^{v \mapsto vA} & \mathbb{Z}^{K} \ar[r]^{v \mapsto vA} & \cdots .
}
$$

We can therefore write
$$
K_0(\SSig) \cong (\mathbb{Z}^{K} \times \mathbb{N})/\! \! \sim ,
$$
where, for $n\leq m$, $(v,n) \sim (w,m)$ if and only if there exists $k \in \mathbb{N}$ such that $vA^{k+m-n} = wA^k$. We write $[v,n]$ for the equivalence class under $\sim$.

Once again proceeding as in the case of $\HSig$ we can show that $\alpha_{\ast}[v,N] = [vA^2,N+1] = [vA,N]$ and $\alpha_{\ast}^{-1}[v,N] = [v,N+1]$.

We are now ready to describe $K_{\ast}(\CH)$ in the SFT case. Once again, this follows \cite{putnam96}. As \newline
$K_1(\HSig) \cong K_0(SH) = 0$ the six term exact sequence from section \ref{KCHsection} becomes 
$$
\xymatrix{
0  \ar[r] & K_0(\CH)\ar[r]^-{(e_0)_{\ast}} & K_0(\HSig)\ar[d]^{(id-\alpha_{\ast})} \\
0  & K_1(\CH) \ar[l]  & K_0(\HSig)\ar[l]^-{\iota_{\ast}}
}
$$

We see immediately that 
$$
K_0(\CH) \cong ker(id-\alpha_{\ast}),
$$
and
$$
K_1(\CH) \cong coker(id-\alpha_{\ast}).
$$
The following well known Lemma will be of use in describing $K_1(\CH)$.

\begin{lemma}\label{LimitSubgroup}
Let $G$ be an abelian group, $\psi: G \rightarrow G$ an endomorphism, and $H < G$ a $\psi$-invariant subgroup, and consider the following diagram
$$
\xymatrix@C+1pc{
H \ar[r]^-{\psi} \ar[d]^-{\iota} & H \ar[r]^-{\psi} \ar[d]^-{\iota} & \cdots \\
G \ar[r]^-{\psi} & G \ar[r]^-{\psi} & \cdots
}
$$
where the vertical maps $\iota$ are given by inclusion, $\iota(x) = x$.
Then $\lim_{\rightarrow}H < \lim_{\rightarrow}G$ and 
$$ \lim_{\rightarrow}G / \lim_{\rightarrow}H \cong  \lim_{\rightarrow}(G/H)$$
\end{lemma}

As indicated in the introduction, we make the following definition.
\begin{definition}\label{C(A)B(A)}
 Let $A$ be a $K \times K$ integer matrix. We define
 $$
C(A) = \{ X \in M_K(\mathbb{Z}) \ | \ AX=XA \},
$$
and
$$
B(A) = \{ AY-YA \mid   Y \in M_K(\mathbb{Z}) \}.
$$
\end{definition}

The following theorem gives a description of $K_0(\CH)$ and $K_1(\CH)$ in terms of inductive limits.

\begin{theo}\label{KCHMatrices}
Let $(\Sigma, \sigma)$ be a mixing SFT with $K\times K$ adjacency matrix $A$ and $\CH$ the corresponding mapping cylinder.

Then $K_0(\CH)$ is isomorphic to the inductive limit of the following system
$$
\xymatrix@C+1pc{
C(A) \ar[r]^-{X \mapsto AXA}  & C(A) \ar[r]^-{X \mapsto AXA} & C(A) \ar[r]  & \cdots 
}
$$
Similarly, $K_1(\CH)$ is isomorphic to the inductive limit of the following system
$$
\xymatrix@C+1pc{
M_K(\mathbb{Z})/B(A) \ar[r]^-{X \mapsto AXA}  & M_K(\mathbb{Z})/B(A) \ar[r]^-{X \mapsto AXA} & M_K(\mathbb{Z})/B(A) \ar[r]  & \cdots 
}
$$
\end{theo}

\paragraph{Proof:}
Consider the following diagram:
$$
\xymatrix@C+1pc{
C(A) \ar[r]^-{X \mapsto AXA} \ar[d]^-{\iota} & C(A) \ar[r] \ar[d]^-{\iota} & \cdots \\
M_n({\mathbb Z}) \ar[r]^-{X \mapsto AXA} & M_n({\mathbb Z}) \ar[r] & \cdots
}
$$
Where the vertical maps are given by $\iota(X) = X$.  The diagram clearly commutes, so $\iota$ induces a well defined group homomorphism $\iota : \lim C(A) \rightarrow K_0(H)$.  It is also clear that $im(\iota) \subset ker(id-\alpha_{\ast})$.  We are left to show that $\iota$ is injective and that the image of $\iota$ is $ker(id-\alpha_{\ast}) \subset K_0(\HSig)$.  

We start by showing injectivity.  Let $[X,m] \in \lim C(A)$ be such that $\iota([X,m]) = [0,m+k] \in K_0(\HSig)$ for some $k$, then $A^{k+l} X A^{k+l} = A^l0A^l$ for some $l$, i.e., $A^jXA^j = 0$ for some $j \geq n$.  But then $[X,m]=[0,j+m] \in \lim C(A)$, and $\iota$ is injective.  

Now suppose $[X,k] \in ker(id-\alpha_{\ast}) \subset K_0(H)$.
Then $(id-\alpha{\ast})[X,k] = [0,m]  \in K_0(H)$ for some $m>k$, i.e., $[AXA-XA^2,k+1] = [0,k]$ or $A^{m-k+l}XA^{m-k+l} - A^{m-k+l-1}XA^{m-k+l+1} = 0$ for some $l$.  Letting $j=m-k+l$ we get $A^{j}XA^{j} = A^{j-1}XA^{j+1}$ or after multiplying on the left by $A$, $A^{j+1}XA^{j} = A^{j}XA^{j+1}$. So $Y = A^{j}XA^{j} \in C(A)$ and $\iota([Y,j+k]) = [X,k]$.  So $\iota: \lim C(A) \rightarrow ker(id - \alpha_{\ast})$ is an isomorphism.

We now show that $im(id-\alpha_{\ast})$ is isomorphic to the inductive limit of
$$
\xymatrix@C+1pc{
B(A) \ar[r]^-{X \mapsto AXA} & B(A) \ar[r]^-{X \mapsto AXA}  & B(A) \ar[r]^-{X \mapsto AXA} & \cdots
}
$$
Consider the diagram
$$
\xymatrix@C+1pc{
B(A) \ar[r]^-{X \mapsto AXA} \ar[d]^-{\iota} & B(A) \ar[r] \ar[d]^-{\iota} & \cdots \\
M_K({\mathbb Z}) \ar[r]^-{X \mapsto AXA} & M_K({\mathbb Z}) \ar[r] & \cdots
}
$$

where the vertical maps, $\iota$, are given by inclusion $\iota(X) = (X)$.  Clearly the above diagram commutes, so $\iota$ extends to a well defined map on the inductive limit groups:
$$\iota : \lim_{\rightarrow} B(A) \longrightarrow \lim_{\rightarrow} M_K({\mathbb Z}) \quad (\cong K_0(H)).$$
We first show that $im(\iota) \subset im(id-\alpha_{\ast})$. 
Suppose $[X,k] \in (B(A) \times {\mathbb N})/\sim $, then there exists $Y$ such that $X = YA - AY$, so 
$$\iota[X,k] = \iota[YA-AY,k] = [YA-AY,k] \in (M_{K}({\mathbb Z}) \times {\mathbb N})/\sim .$$ 
Now 
$$
[YA-AY,k] \! = [A(YA-AY)A,k+1] \! = [(AY)A^2-A(AY)A,k+1] \! = (id-\alpha_{\ast})[-YA,k],
$$ 
so 
$$\iota(\lim_{\rightarrow} B(A)) \subset im(id-\alpha_{\ast}).$$

A straightforward argument now shows that $\iota$ is $1-1$ and onto, so
$$\iota : \lim_{\rightarrow} B(A) \rightarrow im(id-\alpha_{\ast}) $$
is an isomorphism.
Now using lemma \ref{LimitSubgroup} we see immediately that $K_1(\CH) \cong coker(id-\alpha_{\ast}) \cong K_0(\HSig)/im(id-\alpha_{\ast})$ is the limit of the inductive system
$$
\xymatrix@C+1pc{
M_K(\mathbb{Z})/B(A) \ar[r]^-{X \mapsto AXA}  & M_K(\mathbb{Z})/B(A) \ar[r]^-{X \mapsto AXA} & M_K(\mathbb{Z})/B(A) \ar[r]  & \cdots 
}
$$
\qed

\subsection{The Product on $K_{\ast}(\CH)$}\label{RingSection}

 We wish to describe the ring structure on $K_{\ast}(\CH)$ in terms of the matrix characterization of $K_{\ast}(\CH)$ obtained in section \ref{K(CH)}.
We start by determining the ring structure on the subring $K_0(\CH)$.

We begin with a couple of observations.
\begin{enumerate}
\item The real parameter $t$ in Prop. \ref{Kprod} can clearly be replaced with the integer parameter $n$.  This will be helpful as in the SFT case $K_0(\CH) \subset K_0(\HSig)$, and on $K_0(\HSig)$ $\alpha_{\ast}^n$ is defined, but $(\alpha_t)_{\ast}$ is not.
\item In the case of a SFT, we can use a slightly simpler, though equivalent, definition for the product on $K_0(\CH)$. See Prop. \ref{SFTProd}.
\end{enumerate}

Consider the following $\ast$-subalgebra of $\HSig$.
$$
\mathcal{H}(\Sigma,\sigma) = span\{e_{N, v_i, v_j}(\xi, \eta) \}.
$$

Notice that $\mathcal{H}(\Sigma,\sigma)$ is dense in $\HSig$ and that, for each $p \in M_{\infty}(\HSig)$,
there exists \newline 
$q \in M_{\infty}(\mathcal{H}(\Sigma,\sigma))$ such that $[p]_0 = [q]_0$. 
We also consider the following $\ast$-subalgebra of $\CH$, $\mathcal{CH} = \{ f\in \CH \ | \ f(0) \in \mathcal{H} \}$. We again notice that $\mathcal{CH}$ is dense in $\CH$ and, for each $p \in M_{\infty}(\CH)$, there exists $q \in M_{\infty}(\mathcal{CH})$ such that $[p]_0 = [q]_0$.

The following lemma will be useful in computing the product on $K_0(\CH)$.

\begin{lemma}\label{SFTProd1}
Let $e_{N, v_i, v_j}(\xi, \xi)$, $e_{N, v_l, v_k}(\eta, \eta) \in C_c(\GhSig)$ be as in section \ref{C*SFT}. For $n \geq N$ we have
$$
\alpha^n(e_{N, v_i, v_j}(\xi, \xi)) \alpha^{-n}(e_{N v_l v_k}(\eta, \eta)) = \sum_{\xi_1 \in S}e_{N+n,v_i,v_l}(\xi \xi_1 \eta, \xi \xi_1 \eta),
$$
where 
$$
S = \{ \xi_1 \ | \ |\xi_1| = 2n-2N, \ i(\xi_1)=t(\xi), \ t(\xi_1) = i(\eta). \}
$$
\end{lemma}

\paragraph{Proof:}
For notational convenience, let $a = e_{N v_i v_j}(\xi, \xi)$, $b = e_{N v_l v_k}(\eta, \eta)$.  Then 
$$
 \alpha^n(a) = \sum_{v \in V(G)} \sum_{\xi' \in \Xi_{2n,v_j,v}}e_{N+n,v_i,v}(\xi \xi', \xi \xi')
$$
and
$$ 
\alpha^{-n}(b) =  \sum_{v \in V(G)} \sum_{\eta' \in \Xi_{2n,v,v_k}}e_{N+n,v,v_k}(\eta' \eta, \eta' \eta)
$$
so
$$
\alpha^n(a) \alpha^{-n}(b) = \sum_{v \in V(G)} \sum_{\xi'} \sum_{\tilde{v} \in V(G)} \sum_{\eta'}e_{N+n,v_i,v}(\xi \xi', \xi \xi')e_{N+n,v,v_k}(\eta' \eta, \eta' \eta).
$$
Where the sum is over $\xi'$ such that $i{\xi'} = t{\xi}$ and $\eta'$ such that $t(\eta')=i(\eta)$. Furthermore, each summand is $0$ unless $\xi\xi' = \eta'\eta$. 
Write $\xi'=\xi_1 \xi_2$, $\eta'=\eta_1 \eta_2$ where $|\xi_2|=|\eta|=|\xi|=|\eta_1|=2N$ and $|\xi_1|=|\eta_2|=2n-2N$. Now $\xi \xi_1 \xi_2 = \eta_1 \eta_2 \eta$ implies $\xi = \eta_1$, $\xi_2=\eta$, and $\xi_1 = \eta_2$, which in turn imply $v_i = \tilde{v}$, $v=v_l$, $i(\xi_1)=i(\eta_2)=t(\xi)=v_j$, and $t(\xi_1)=t(\eta_2)=i(\eta)=v_k.$ So the sum becomes
$$
\alpha^n(a) \alpha^{-n}(b) = \sum_{\xi_1\in S}e_{N+n,v_i,v_l}(\xi \xi_1 \eta),
$$
where
$$
S = \{ \xi_1 \ | \ |\xi_1| = 2n-2N, \ i(\xi_1)=t(\xi), \ t(\xi_1) = i(\eta). \}
$$
\qed

\begin{prop}\label{SFTProd}
Let $(\Sigma, \sigma)$ be an irreducible SFT, $p,q \in P_{\infty}(\mathcal{CH})$. By identifying $K_0(\CH)$ with a subgroup of $K_0(\HSig)$, the product can be written 
$$
[p(0)]_0[q(0)]_0 = \lim_{n \rightarrow \infty}[\alpha^n(p(0))\alpha^{-n}(q(0))]_0.
$$
\end{prop}

\paragraph{Proof:}
From section \ref{K(CH)} we know that $K_0(\CH) \cong ker(id-\alpha_{\ast}) \subset K_0(\HSig)$ where the isomorphism is given by evaluation at $0$. So we have
\begin{eqnarray*}
[p(0)]_0[q(0)]_0 &=& \lim_{t \rightarrow \infty} \left[\chi_{(1/2,\infty)}(p \times q)_t(0) \right]_0 \\
&=& \lim_{t \rightarrow \infty} \left[\chi_{(1/2,\infty)}(p(0) \times q(0))_t \right]_0 \\
&=& \lim_{n \rightarrow \infty} \left[\chi_{(1/2,\infty)}(p(0) \times q(0))_n \right]_0.
\end{eqnarray*}
We now look more closely at 
$$
(p(0) \times q(0))_n = \frac{\alpha^n(p(0))\alpha^{-n}(q(0))+ \alpha^{-n}(q(0))\alpha^n(p(0))}{2}.
$$
Since $p(0), q(0) \in M_{\infty}(\mathcal{H}(\Sigma,\sigma))$, it will suffice to consider the product $\alpha^n(a)\alpha^{-n}(b)$, where $a = e_{N v_i v_j}(\xi, \xi)$  and $b=e_{N v_l v_k}(\eta, \eta)$.
Fix $n \geq N$, then from Lemma \ref{SFTProd1}
$$
\alpha^n(a) \alpha^{-n}(b) = \sum_{\xi_1\in S}e_{N+n,v_i,v_l}(\xi \xi_1 \eta).
$$
A similar computation yields
$$
\alpha^{-n}(b)\alpha^n(a) = \sum_{\eta_2 \in S}e_{N+n,v_i,v_l}(\xi \eta_2 \eta),
$$
so these two sums are equal. This implies 
$$
\alpha^n(a) \alpha^{-n}(b) = \alpha^{-n}(b)\alpha^n(a) 
$$
is a projection.
Now for $p$, $q$ as above, there exists $N$ such that for $n>N$
$$
\alpha^n(p(0)) \alpha^{-n}(q(0)) = \alpha^n(q(0)) \alpha^{-n}(p(0)),
$$
and this is a projection. Thus we have
$$
[p(0)]_0[q(0)]_0 =  \lim_{n \rightarrow \infty}[\alpha^n(p(0))\alpha^{-n}(q(0))]_0.
$$
\qed

\begin{remark}
The preceding proposition gives the product on positive elements in $K_0$.  The product extends to all of $K_0$ by linearity and it is clear that it is positive.
\end{remark}

Under the isomorphism of Theorem \ref{KCHMatrices}, every element of $K_0(\CH)$ is equal to some $[X,N] \in \lim C(A)$. As each such $X$ is a linear combination of matrices of the form $e_{ij}$, we start with two matrices of this form, their corresponding projections in $\HSig$ and multiply according to Prop. \ref{SFTProd}.
\begin{remark}
As $[e_{ij},N] \in K_0(\HSig)$ need not be an element of the subgroup $K_0(\CH)$, the formula we derive for the product of two such elements will not be well defined in general.  I.e., the element of $K_0(\HSig)$, $[\alpha^n(a)\alpha^{-n}(b)]_0$, in the following Lemma depends (in general) on the integer $n$ and thus $\lim_{n \rightarrow \infty}[\alpha^n(a)\alpha^{-n}(b)]_0$ need not exist.  However, if we apply the formula to linear combinations of such elements, $[X,N]$, $[Y,M]$ which are in $K_0(\CH)$ the product is well defined.
\end{remark}

\begin{lemma}\label{HProd}
Let $a, b \ \in P_{\infty}(\HSig)$.  If $[a]_0 = [X,N]$, $[b]_0 = [Y,N]$, then for $n \geq N$
$$
[\alpha^n(a)\alpha^{-n}(b)]_0 = [XA^{2n-2N}Y, N+n].
$$
\end{lemma}

\paragraph{Proof:}
For $a, b \ \in P_{\infty}(\HSig)$, we can find $\bar{a}, \bar{b} \ \in P_{\infty}(\mathcal{H}(\Sigma,\sigma))$ such that $[a]_0 = [\bar{a}]_0$, $[b]_0 = [\bar{b}]_0$.
It therefore suffices to prove the result for rank one projections $\bar{a} = e_{N v_i v_j}(\xi, \xi)$  and $\bar{b}=e_{N v_l v_k}(\eta, \eta)$. Then $[\bar{a}]_0 = [e_{ij},N], [\bar{b}]_0 = [e_{kl},N] \in K_0(\HSig)$.
Fix $n \geq N$, then using Lemmas \ref{SFTProd1} and \ref{SFTProd} we see that
$$
\alpha^n(\bar{a}) \alpha^{-n}(\bar{b}) = \sum_{\xi_1}e_{N+n,v_i,v_l}(\xi \xi_1 \eta)
$$
is a projection.
The number of summands is thus the number of paths $\xi_1$ of length $2n-2N$ from $v_j$ to $v_k$, i.e., $A^{2n-2N}_{jk}$.
Noticing that $A^{2n-2N}_{jk} e_{il} = e_{ij}A^{2n-2N}e_{kl}$, we have 
$$
[\alpha^n(\bar{a}) \alpha^{-n}(\bar{b})]_0 = [e_{ij}A^{2n-2N}e_{kl}, N+n].
$$
\qed

We are now ready to prove Theorem \ref{K0XK0}.

\paragraph{Proof of Theorem \ref{K0XK0}:}
Let $p, \ q \in P_{\infty}(\mathcal{CH})$, with $[p(0)]_0 = [X,N]$, $[q(0)]_0 = [Z,N]$ (under the isomorphism in Prop. \ref{KCHMatrices}). Note here that we have chosen representatives of each equivalence class such that the integer $N$ is equal. From Prop. \ref{SFTProd} we know the product is given by
$$
[p(0)]_0[q(0)]_0 = \lim_{n \rightarrow +\infty}[\alpha^n(p(0))\alpha^{-n}(q(0))]_0.
$$
However, we know that for $n \geq N$ the sequence $[\alpha^n(p(0))\alpha^{-n}(q(0))]_0$ is constant, so we can write
$$
[p]_0[q]_0 = \lim_{n \rightarrow +\infty}[\alpha^n(p(0))\alpha^{-n}(q(0))]_0 = [\alpha^{N}(p(0))\alpha^{-N}(q(0))]_0
$$
By Lemma \ref{HProd} this is
$$
[X,N]\ast[Z,N] = \lim_{n \rightarrow +\infty}[XA^{2n-2N}Z,N+n] = [XZ,2N].
$$
Now suppose $[Y,M] \in K_0(\CH)$.  Let $M_1 = \max\{N,M\}$.
\begin{eqnarray*}
[X,N]\ast[Y,M] &=& [A^{M_1-N}XA^{M_1-N},M_1]\ast[A^{M_1-M}YA^{M_1-M},M_1] \\
&=& [A^{M_1-N}XA^{M_1-N}A^{M_1-M}YA^{M_1-M}, 2M_1 ] \\
&=& [A^{M_1-N}A^{M_1-M}XYA^{M_1-M}A^{M_1-N}  , N+M + |M-N| ] \\
&=& [A^{|M-N|}XYA^{|M-N|}, N+M+|M-N|] \\
&=& [XY,N+M]
\end{eqnarray*}
\qed

\begin{cor}
The multiplicative identity $[I,0]$ in the ordered group $K_0(\CH)$ is  an order unit. 
\end{cor}

\paragraph{Proof:}
That $[I,0]$ is the multiplicative identity is immediate from Theorem \ref{K0XK0}.
As $(\Sigma, \sigma)$ is mixing, $A$ is primitive, so there exists $K$ such that $A^{2k}$ is strictly positive for all $k>K$.  Let $[X,k] \in K_0(\CH)$, wolog $k>K$. Now there exists an integer $l$ such that $lA^{2k} > X$ entry-wise. In other words $l[I,0] = l[A^{2k},k] > [X,k]$.
\qed

Before going on to consider the other cases of the product, we consider an example. Consider the SFT, $(\Sigma,\sigma)$ with adjacency matrix
$$
A = \left[\begin{array}{ccc}
2 & 1 & 1 \\
1 & 2 & 1 \\
1 & 1 & 2 \end{array} \right]. 
$$

It is straightforward to check that  $C(A)$ consists of integral combinations of the matrices:

\begin{align*}
X_1&=\left[\begin{array}{ccc}
1 & -1 & 0 \\
-1 & 1 & 0 \\
0 & 0 & 0 
\end{array}\right]
& X_2&=\left[\begin{array}{ccc}
0 & 1 & -1 \\
0 & -1 & 1 \\
0 & 0 & 0 
\end{array}\right]
& X_3&=\left[\begin{array}{ccc}
0 & 0 & 0 \\
1 & -1 & 0 \\
-1 & 1 & 0 
\end{array}\right] \\
X_4&=\left[\begin{array}{ccc}
0 & 0 & 0 \\
0 & 1 & -1 \\
0 & -1 & 1 
\end{array}\right]
& X_5&=\left[\begin{array}{ccc}
1 & 0 & 0 \\
0 & 1 & 0 \\
0 & 0 & 1 
\end{array}\right] = I
\end{align*}

The trace map $\tau^{CH}$ maps $K_{0}(\CH)$ onto $\mathbb Z[1/2]$. The kernel is a free abelian group and is generated by $[X_{i}, 0], 1 \leq i \leq 4$.
It is easy to check that this ideal is non-commutative.

We now calculate the product of an element of $K_0(\CH)$ with an element of $K_1(\CH)$.  First recall the 6-term exact sequence:
$$
\xymatrix{
0 \ar[r] & K_0(\CH)\ar[r]^-{(e_0)_{\ast}} & K_0(\HSig)\ar[dr]^{id-\alpha_{\ast}} \\
0 & K_1(\CH) \ar[l] & K_1(SH) \ar[l]_-{\iota_{\ast}} & K_0(\HSig)\ar[l]_-{\cong}
}
$$

Hence every element of $K_1(\CH)$ is the image under $\iota_{\ast}$ of some element in $K_1(SH) \cong K_0(\HSig)$. Also recall that $K_0(H)$ is generated by rank one projections in $P_1(\HSig)$.  We will proceed as follows.  Starting with a $p \in P_1(\HSig)$, and $q \in P_{m}(\CH)$. We find $u_p \in U_1(\widetilde{SH})$ (corresponding to the isomorphism $K_1(SH) \cong K_0(\HSig)$). We then find $\iota_{\ast}(u_p) \in U_1(\widetilde{\CH})$ and $\tilde{p} \in P_2(\widetilde{S\CH})$ (corresponding to the isomorphism $K_1(\CH) \cong K_0(S\CH)$).  We then simply embed $\tilde{p}$ and $q$ in $P_{\infty}(C(S^1,\CH))$ and multiply according to Lemma \ref{SFTProd2}.

The following Lemmas are standard results in $K$-theory, see for example Theorems 10.1.3 and 11.1.2 in \cite{rordam}.

\begin{lemma}\label{K0ToK1}
Let $p \in P_m(H)$, then under the isomorphism $K_0(H) \cong K_1(SH)$, $p \mapsto u_p \in U_m(\widetilde{SH})$ where
$$
u_p(s) = e^{2\pi is}p + (1-p), \ \textrm{for $0\leq s \leq 1$}
$$
\end{lemma}

\begin{lemma}\label{K1toK0}
Let $u_p \in U_m(\widetilde{\CH})$, then under the isomorphism $K_1(\CH) \cong K_0(\widetilde{S\CH})$, $u_p \mapsto \tilde{p} \in P_{2m}(\widetilde{S\CH})$ where
$$
\tilde{p} = v_p \left [ 
\begin{array}{cc}
I_m & 0_m \\
0_m & 0_m 
\end{array} \right ] v_p^{\ast}
$$
and 
$$
v_p(t) = \overline{R_t}\left [ 
\begin{array}{cc}
u_p & 0 \\
0 & I_m 
\end{array} \right ]
\overline{R_t}^{\ast} \left [ 
\begin{array}{cc}
I_m & 0 \\
0 & u_p^{\ast} 
\end{array} \right ] $$
where 
$$
\overline{R_t} = \left[ \begin{array}{cc}
\cos(\frac{\pi}{2} t)I_m & -\sin(\frac{\pi}{2} t)I_m \\
\sin(\frac{\pi}{2} t)I_m & \cos(\frac{\pi}{2} t)I_m 
\end{array}
\right ],
$$
so that
$$
v_p(0) = \left [ 
\begin{array}{cc}
u_p & 0 \\
0 & u_p^{\ast} 
\end{array} \right ], \quad v_p(1) = \left [ 
\begin{array}{cc}
I_m & 0 \\
0 & I_m 
\end{array} \right ].
$$
\end{lemma}

We state and prove one more Lemma before we compute the product of an element of $K_0(\CH)$ with an element of $K_1(\CH)$.

\begin{lemma}\label{HXCH}
For $p \in P_1(\mathcal{H}(\Sigma,\sigma))$, $q \in P_{\infty}(\mathcal{CH})$,
there exists $N \in \mathbb{N}$ such that, for $n, \ m \geq N$, the matrix $(p \times q)_n$ with $(i,j)$ entry given by $\alpha^n(p)\alpha^{-n}(q_{ij}(0))$ is in $P_{\infty}(\HSig)$ and
$$
[(p\times q)_n]_0 - [(p \times q)_m]_0 \in Im(id - \alpha_{\ast}).
$$
\end{lemma}
\paragraph{Proof:}
The existence of $N$ such that $(p \times q)_n \in P_{\infty}(H)$ for all $n \geq N$ follows immediately from \ref{SFTProd}.  Now fix $m > n \geq N$ and consider
\begin{eqnarray*}
\alpha_{\ast}[\alpha^n(p)\alpha^{-n}(q_{ij}(0))]_0 &=& [\alpha^{n+1}(p)\alpha^{-n+1}(q_{ij}(0))]_0 \\
&=& [\alpha^{n+1}(p)\alpha^{-n+1}(q_{ij}(-2))]_0 \ \ \textrm{by homotopy invariance} \\
&=& [\alpha^{n+1}(p)\alpha^{-n-1}(q_{ij}(0))]_0.
\end{eqnarray*}
So
$$
(id -\alpha_{\ast})[\alpha^n(p)\alpha^{-n}(q_{ij}(0))]_0 = [\alpha^n(p)\alpha^{-n}(q_{ij}(0))]_0 - [\alpha^{n+1}(p)\alpha^{-n-1}(q_{ij}(0))]_0,
$$
and hence, by induction
$$
[\alpha^n(p)\alpha^{-n}(q_{ij}(0))]_0 - [\alpha^{m}(p)\alpha^{-m}(q_{ij}(0))]_0 \in Im(id-\alpha_{\ast}).
$$
\qed
\begin{remark}
The same result clearly holds for the product $(q \times p)_n$ with entries $\alpha^n(q_{ij}(0))\alpha^{-n}(p)$.
\end{remark}

\begin{definition}
For $p \in P_1(\mathcal{H}(\Sigma,\sigma))$, $q \in P_{\infty}(\mathcal{CH})$, $n \geq N$ define $ (p \times q)_n$ as in Lemma \ref{HXCH}.
With a slight abuse of notation, we will often drop the $n$ and write $(p \times q)$. 
\end{definition}

The proof of Theorem \ref{K0XK1}
is quite long, so we will break the proof down into a series of Lemmas.

\begin{lemma}\label{UpQ}
Let $p \in P_1(\mathcal{H}(\Sigma,\sigma))$, $q \in P_{m}(\mathcal{CH})$ for some $m$. Let $u_p \in U_1(\widetilde{\CH})$ be
$$
u_p(s) = e^{2\pi is}\alpha^n(p) + (1-\alpha^n(p)) \quad \textrm{for $n \leq s \leq n+1$},
$$
Then for $n$ large enough (as in Lemma \ref{HXCH}), the matrix $(u_p \times q)_n$ is given by 
$$
((u_p \times q)_n)_{ij} = \alpha_n(u_p)\alpha_{-n}(q_{ij}) = (u_{p\times q})_{ij} - (I_m - \alpha_{-n}(q_{ij})).
$$
in other words
$$
(u_p \times q)_n = u_{p\times q} - (I_m - \alpha_{-n}(q)) = u_{p\times q}\alpha_{-n}(q).
$$
\end{lemma}
\paragraph{Proof:}
\begin{eqnarray*}
\left(\alpha_n(u_p)\alpha_{-n}(q_{ij})\right)(s) 
 \! & = &\!  u_p(s+n)q_{ij}(s-n) \\
 \! & = &\!  \alpha_{s}\left((u_p(n))(q_{ij}(-n))\right) \\
 \! & = &\!  \alpha_{s}\left(\big(e^{2\pi in}\alpha^n(p) + (1-\alpha^n(p))\big)\alpha^{-n}(q_{ij}(0)) \right) \\
 \! & = &\!  \alpha_{s}\left(e^{2\pi in}\alpha^n(p)\alpha^{-n}(q_{ij}(0)) -\alpha^n(p)\alpha^{-n}(q_{ij}(0)) + \alpha^{-n}(q_{ij}(0)) \right) \\
 \! & = &\!  \alpha_{s}\left(e^{2\pi in}((p \times q)_n)_{ij} - ((p \times q)_n)_{ij} + \alpha^{-n}(q_{ij}(0)) \right) \\ 
 \! & = &\!  \alpha_{s}\left(e^{2\pi in}(p \times q)_n + (I_m - (p \times q)_n) - (I_m - \alpha^{-n}(q(0))) \right)_{ij} \\
 \! & = &\!  \alpha_s \left( u_{(p \times q)}(0) - (I_m - \alpha^{-n}(q(0))) \right)_{ij} \\
 \! & = &\!  \left( u_{(p \times q)}(s) - (I_m - \alpha^{-n}(q(s))) \right)_{ij}
 \end{eqnarray*}
so 
$$
(u_p \times q)_n = u_{p\times q} - (I_m - \alpha_{-n}(q)).
$$
Now
\begin{eqnarray*}
\alpha_n(u_p)\alpha_{-n}(q) &=& \alpha_n(u_p)\alpha_{-n}(q)^2 \\
&=& \left( u_{p\times q} - (I_m - \alpha_{-n}(q))\right)\alpha_{-n}(q) \\
&=& u_{p\times q}\alpha_{-n}(q).
\end{eqnarray*}
\qed

\begin{lemma}
Let $p$ and $q$ be as in Lemma \ref{UpQ}, $u_{p\times q}$ and $v_{p\times q}$ as in Lemmas \ref{K0ToK1} and \ref{K1toK0}. Then
$$
u_{p\times q}(I_m - \alpha_{-n}(q)) = I_m - \alpha_{-n}(q)
$$
and
$$
v_{p\times q}\left[\begin{array}{cc}
I_m - \alpha_{-n}(q) & 0 \\
0 & I_m - \alpha_{-n}(q)
\end{array}\right] = \left[\begin{array}{cc}
I_m - \alpha_{-n}(q) & 0 \\
0 & I_m - \alpha_{-n}(q)
\end{array}\right]
$$
\end{lemma}

\paragraph{Proof:}
Straightforward calculations from the definitions of $u_{p\times q}$, $v_{p\times q}$ and the result of Lemma \ref{UpQ}. \qed

The following Lemma allows us to use a slightly simplified form of the product, similar to Prop. \ref{SFTProd}.

\begin{lemma}\label{SFTProd2}
Let $p \in P_1(\mathcal{H}(\Sigma,\sigma))$, $q \in P_{m}(\mathcal{CH})$. Let $\tilde{p} \in P_2(C(S^1,\CH))$, then there exists $T \in \mathbb{R}$ such that for $t>T$
$$
\alpha_t(\tilde{p})_{(ij)}\alpha_{-t}(q_{kl}) = \alpha_{-t}(q_{kl})\alpha_t(\tilde{p}_{(ij)}).
$$ 
Moreover, if we define $(\tilde{p} \times q)_t$ componentwise by  
$$
((\tilde{p} \times q)_t)_{(ik)(jl)} = \alpha_t(\tilde{p}_{(ij)})\alpha_{-t}(q_{kl})
$$
then for $t>T$ $(\tilde{p} \times q)_t$ is a projection in $P_{2m}(C(S^1,\CH))$ and
$$
[\tilde{p}]_0[q]_0 = \lim_{t \rightarrow \infty}[(\tilde{p} \times q)_t]_0
$$
\end{lemma}

\paragraph{Proof:}

The product $(\tilde{p} \times q)_n$ defined in the statement of the Lemma is the same as regular matrix multiplication between the following two $2m \times 2m$ matrices.
$$
(\tilde{p} \times q)_n = \left [ \begin{array}{cc}
\alpha_n{\tilde{p}_{11}}I_m & \alpha_n{\tilde{p}_{12}}I_m \\
\alpha_n{\tilde{p}_{21}}I_m & \alpha_n{\tilde{p}_{22}}I_m
\end{array}\right ]
\left [\begin{array}{cc}
\alpha_{-n}(q) & 0 \\
0 & \alpha_{-n}(q)
\end{array}\right ].
$$
For the $2\times 2$ matrix $Y$ we will denote by $\bar{Y}$ the $2m \times 2m$ with $m \times m$ block entries $Y_{ij}I_m$.  Also notice that $\overline{XY} = \bar{X}\bar{Y}$. We now have
\begin{eqnarray*}
(\tilde{p} \times q)_n &=&  \alpha_n\Big(\bar{v_p}\left [ 
\begin{array}{cc}
I_m & 0 \\
0 & 0 
\end{array} \right ] \bar{v_p^{\ast}}\Big)\alpha_{-n}\left [ 
\begin{array}{cc}
q & 0 \\
0 & q 
\end{array} \right ] \\
& =& 
\overline{\alpha_n(v_p)}\left [ 
\begin{array}{cc}
I_m & 0 \\
0 & 0 
\end{array} \right ] \overline{\alpha_n(v_p^{\ast})}\left [ 
\begin{array}{cc}
\alpha_{-n}(q) & 0 \\
0 & \alpha_{-n}(q) 
\end{array} \right ] .
\end{eqnarray*}

From the proof of Lemma \ref{UpQ} and Prop. \ref{SFTProd} we see that
$$
\alpha_n(u_p^{\ast})\alpha_{-n}(q) = \alpha_{-n}(q)\alpha_n(u_p^{\ast}).
$$ 

From the definition of $(v_p^{\ast})$, and the fact that $R_t$ is a matrix of scalars, we can show that 
$$
\left [ 
\begin{array}{cc}
\alpha_{-n}(q) & 0 \\
0 & \alpha_{-n}(q) 
\end{array} \right ]
$$
commutes with both $\overline{\alpha_n(v_p)}$ and $\overline{\alpha_n(v_p^{\ast})}$.

So
\begin{eqnarray*}
(\tilde{p} \times q)_n &=& 
 \overline{\alpha_n(v_p)}\left [ 
\begin{array}{cc}
I_m & 0 \\
0 & 0 
\end{array} \right ] \overline{\alpha_n(v_p^{\ast})}\left [ 
\begin{array}{cc}
\alpha_{-n}(q) & 0 \\
0 & \alpha_{-n}(q) 
\end{array} \right ] \\
&=& \left [ 
\begin{array}{cc}
\alpha_{-n}(q) & 0 \\
0 & \alpha_{-n}(q) 
\end{array} \right ] 
 \overline{\alpha_n(v_p)}\left [ 
\begin{array}{cc}
I_m & 0 \\
0 & 0 
\end{array} \right ] \overline{\alpha_n(v_p^{\ast})}.
\end{eqnarray*}
Hence for sufficiently large $n$ $(\tilde{p} \times q)_n \in P_{2m}(C(S^1,\CH))$ and 
$$
[\tilde{p}]_0[q]_0 = \lim_{n \rightarrow \infty}[(\tilde{p} \times q)_n]_0
$$
\qed

The following Lemma contains the bulk of the calculation involved in proving  \ref{K0XK1}.

\begin{lemma}\label{K0XK1:1}
Let $p \in P_1(\mathcal{H}(\Sigma,\sigma))$, $q \in P_m(\mathcal{CH})$. Let $u_p \in U_1(\widetilde{\CH})$ be
$$
u_p(s) = e^{2\pi is}\alpha^n(p) + (1-\alpha^n(p)) \quad \textrm{for $n \leq s \leq n+1$}.
$$
Let $\tilde{p} \in P_2(\widetilde{S\CH})$ be as in Lemma \ref{K1toK0}. Then 
$$
[\tilde{p}]_0[q]_0 = [ \widetilde{p\times q} ]_0 - [I_m - q]_0
$$
\end{lemma}

\paragraph{Proof:}
As $\tilde{p} \in P_2(\widetilde{C(S^1,CH)})$, $q \in P_m(\tilde{C(S^1,CH)})$, by Lemma \ref{SFTProd2} the product is (for sufficiently large $n$) 
$$
[\tilde{p}]_0[q]_0 = [(\tilde{p} \times q)_n]_0.
$$
Where $(\tilde{p} \times q)_n \in P_{2m}(\widetilde{C(S^1,CH)})$ is as defined in Lemma \ref{SFTProd2}. 
Now as in the proof of Lemma \ref{SFTProd2}, 
$$
(\tilde{p} \times q)_n =  \overline{\alpha_n(v_p)}\left [ 
\begin{array}{cc}
I_m & 0 \\
0 & 0 
\end{array} \right ] \overline{\alpha_n(v_p^{\ast})}\left [ 
\begin{array}{cc}
\alpha_{-n}(q) & 0 \\
0 & \alpha_{-n}(q) 
\end{array} \right ],
$$
and
\begin{eqnarray*}
\overline{\alpha_n(v_p^{\ast})}\left [ 
\begin{array}{cc}
\alpha_{-n}(q) & 0 \\
0 & \alpha_{-n}(q) 
\end{array} \right ] &=& \left [ 
\begin{array}{cc}
I_m & 0 \\
0 & \alpha_n(u_p)I_m 
\end{array} \right ]
\bar{R_t} \left [ 
\begin{array}{cc}
\alpha_n(u_p^{\ast})\alpha_{-n}(q) & 0 \\
0 & \alpha_{-n}(q) 
\end{array} \right ] \bar{R_t^{\ast}} \\
&=& \left [ 
\begin{array}{cc}
I_m & 0 \\
0 & \alpha_n(u_p)I_m 
\end{array} \right ]
\bar{R_t} \left [ 
\begin{array}{cc}
u^{\ast}_{p\times q}\alpha_{-n}(q) & 0 \\
0 & \alpha_{-n}(q) 
\end{array} \right ] \bar{R_t^{\ast}} \\
&=& \left [ 
\begin{array}{cc}
\alpha_{-n}(q) & 0 \\
0 & u_{p\times q}\alpha_{-n}(q)
\end{array} \right ]
\bar{R_t} \left [ 
\begin{array}{cc}
u^{\ast}_{p\times q} & 0 \\
0 & I_m 
\end{array} \right ] \bar{R_t^{\ast}} \\
&=& \left [ 
\begin{array}{cc}
\alpha_{-n}(q) & 0 \\
0 & \alpha_{-n}(q)
\end{array} \right ]v^{\ast}_{p \times q}.
\end{eqnarray*}

So we have
\begin{eqnarray*}
(\tilde{p} \times q)_n &=&  \overline{\alpha_n(v_p)}\left [ 
\begin{array}{cc}
I_m & 0 \\
0 & 0 
\end{array} \right ] \overline{\alpha_n(v_p^{\ast})}\left [ 
\begin{array}{cc}
\alpha_{-n}(q) & 0 \\
0 & \alpha_{-n}(q) 
\end{array} \right ] \\
&=& \overline{\alpha_n(v_p)}\left [ 
\begin{array}{cc}
I_m & 0 \\
0 & 0 
\end{array} \right ]\left [ 
\begin{array}{cc}
\alpha_{-n}(q) & 0 \\
0 & \alpha_{-n}(q)
\end{array} \right ]v^{\ast}_{p \times q} \\
&=& v_{p\times q} \left [ 
\begin{array}{cc}
\alpha_{-n}(q) & 0 \\
0 & 0
\end{array} \right ]v^{\ast}_{p \times q} \\
&=& v_{p\times q} \left [ 
\begin{array}{cc}
I_m & 0 \\
0 & 0
\end{array} \right ]v^{\ast}_{p \times q} - v_{p\times q} \left [ 
\begin{array}{cc}
I_m - \alpha_{-n}(q) & 0 \\
0 & 0
\end{array} \right ]v^{\ast}_{p \times q} \\
&=&
\widetilde{p\times q} - \left [ 
\begin{array}{cc}
I_m - \alpha_{-n}(q) & 0 \\
0 & 0
\end{array} \right ].
\end{eqnarray*}
So,
$$
[\tilde{p}]_0[q]_0 = [\tilde{p} \times q]_0 = [\widetilde{p \times q}]_0 - [I_m - \alpha_{-n}(q)]_0 = [\widetilde{p \times q}]_0 - [I_m - q]_0
$$
\qed

\noindent
In the above, $[\tilde{p}]_0$ is an element of $K_0(\widetilde{S\CH})$. We wish to consider $K_0(S\CH)$, for which it suffices to consider elements of the form $[\tilde{p}]_0 - [\tilde{p}(0)]_0$.

\begin{lemma}\label{K0XK1:3}
Let $p \! \in \! P_1(\mathcal{H}(\Sigma,\sigma))$, $q \! \in \! P_m(\mathcal{CH})$ so that
$[\tilde{p}]_0 - [\tilde{p}(0)]_0 \! \in \! K_0(S\CH)$, $[q]_0 \in K_0(\CH)$. Then
$$
([\tilde{p}]_0 - [\tilde{p}(0)]_0)[q]_0 = [\widetilde{p\times q}]_0 - [\widetilde{p\times q}(0)]_0.
$$
\end{lemma}

\paragraph{Proof:}
\begin{eqnarray*}
([\tilde{p}]_0 - [\tilde{p}(0)]_0)[q]_0 &=& ([\tilde{p}]_0 - [I_1]_0)[q]_0 \\
&=& [\tilde{p} \times q]_0 - [q]_0 \\
&=& [\widetilde{p \times q}]_0 - [I_m-q]_0 - [q]_0 \ \ \textrm{Lemma \ref{K0XK1:1}} \\
&=& [\widetilde{p\times q}]_0 - [I_m]_0 \\
&=& [\widetilde{p\times q}]_0 - [\widetilde{p\times q}(0)]_0.
\end{eqnarray*}
\qed

\begin{remark}
The map $\phi: K_0(H) \cong K_1(SH) \hookrightarrow K_1(\CH) \cong K_0(S\CH)$ is given by
$$
\phi([p]_0) = [\tilde{p}]_0 - [\tilde{p}(0)]_0.
$$
Where, for $\tilde{p} \in P_{2n}(S\CH)$, $\tilde{p}(0) = I_n$.
So Lemma \ref{K0XK1:3} simply says that to multiply $\phi([p]_0)$ by $[q]_0$, we can first multiply $[p]_0$ and $[q]_0$ in $K_0(H)$, which by Lemma \ref{HXCH} is well defined modulo $Im(id-\alpha_{\ast})$, and then apply $\phi$. As $Im(id-\alpha_{\ast})$ is exactly the kernel of the map $ K_1(SH) \hookrightarrow K_1(\CH)$ everything is well defined.
\end{remark}
\begin{remark}
Completely analogous calculations to those above show that
$$
[q]_0([\tilde{p}]_0 - [\tilde{p}(0)]_0) = [\widetilde{q\times p}]_0 - [\widetilde{q\times p}(0)]_0.
$$
\end{remark}

The proof of Theorem \ref{K0XK1} is now immediate.

Finally, we turn to the proof of Theorem \ref{K1XK1} which states 
that the  product on $K_{1}(\CH) \times K_{1}(\CH)$ to $K_{0}(\CH)$ is trivial.

\paragraph{Proof of Theorem \ref{K1XK1}:}
We know that $SH$ is an ideal in $\CH$ and the map induced by inclusion is surjective on $K_{1}$.
Secondly, it is clear that for integer values of the parameter $t$, the asymptotic homomorphism
on a pair from $SH$ takes values again in  $SH$. Hence, we conclude that the range of the product
must lie in the image of $K_{0}(SH)$ in $K_{0}(\CH)$ under the inclusion map. On the other hand, 
 $K_{0}(SH) =0$. This completes the proof. 
\qed

The proof of Theorem \ref{module} computing the module structures 
is done in much the same way. In fact, many of the ingredients are already done. But we omit the remainder of the computations.

\section{Duality}\label{secDuality}

Let $A$ have minimal polynomial $m_{A}(x) = x^l(x^k + a_{k-1}x^{k-1} + \cdots + a_0)$ with $a_0 \neq 0$ (so $m_{A}$ has degree $k+l$, and $l$ is the multiplicity of $0$ as a root). Further, we let $p_{A}(x) = x^k + a_{k-1}x^{k-1} + \cdots + a_0$, so 
that $m_{A}(x) = p_{A}(x) x^{l}$.

The following Lemma sharpens our description of $R_{A}$.

\begin{lemma}\label{Subring-R}
Let $\mathbb{Z}_{k-1}[x]$ denote the set of integral polynomials with degree less than $k$ (the degree of $p_{A}(x)$).
Then we have 
$$
 R_{A} = \{ [p(A), N] \in K_{0}(\CH) \mid p \in \mathbb{Z}_{k-1}[x], N \geq 0 \}.
$$
Moreover, for $ p \in \mathbb{Z}_{k-1}[x]$ and $ N \geq 0$, $ [p(A), N]  = 0$ if and only if 
$p=0$.
\end{lemma}

\paragraph{Proof:}
The containment $\supset$ follows from the observation that, for any $i \geq 0$ and $N \geq 0$, we have 
$$
[A^{i}, N] = [A^{i+2j}, N+j] = [A, 0]^{i+j-N} \ast [A, 1]^{N+j},
$$
where $j$ is chosen so that $i+j -N > 0$.

For the reverse inclusion, it is clear that any element of $R_A$ must be of the form 
$[p(A), N]$, where $p$ is some integer polynomial and $N \geq 0$. 
Let $j$ be the degree of $p$. Of course, if $j < k$, then we are done,
so suppose that $j \geq k$. Consider $p(x)x^{2l}$ and divide it by $m_{A}(x)x^{l} = p_{A}(x)x^{2l}$ and write the result:
$$
p(x)x^{2l} = p_{A}(x)x^{2l} d(x) + q(x),
$$
where $q(x)$ has degree less than that of $p_{A}(x)x^{2l}$, namely $2l +k$. 
It follows that $p(A)A^{2l} = q(A)$. Also,
$x^{2l}$ clearly divides both $p(x)x^{2l}$  and $ p_{A}(x)x^{2l} d(x)$, so it also divides $q(x)$. We write
$q(x) = q_{0}(x)x^{2l}$. Now $q_{0}(x)$ has degree less than $k$ and 
$$
[p(A), N] = [p(A)A^{2l}, N+l] = [q(A), N+l] = [q_{0}(A)A^{2l}, N+l] = [q_{0}(A),N].
$$
This completes the proof of the first statement.

For the second statement, the ``if'' part is clear.
Now suppose $[c_{k-1}A^{k-1} + \cdots + c_0I,N] = [0,0]$.  In other words, there exists $m$ 
such that $(c_{k-1}A^{k-1} + \cdots + c_0I)A^{2m} = 0$.  But since $l$ is the multiplicity of $0$ as a 
root of the minimal polynomial, it must be true that $(c_{k-1}A^{k-1} + \cdots + c_0I)A^l = 0$.  So we have $c_{k-l-1}A^{k+l-1} + \cdots + c_0A^l = 0$.  
We recall that the minimal polynomial of $A$ has degree $k+l$, so $\{A^{k+l-1},A^{k+l-2}, \ldots, A^l \}$ is a linearly independent set, and we see that 
all of the $c_i$'s equal $0$.
\qed

\begin{lemma}\label{ModuleHom1}
Let $\varphi$ be in $Hom_{R_{A}}(K_0(S),R_{A})$. Then 
 there exists $z \in \mathbb{Z}^{K}$ (considered as a column vector) and $N \in \mathbb{N}$ such that, for each $[v,n] \in K_0(S)$,

\begin{equation*}
\begin{split}
\varphi[v,n] &= [vA^{n}zA^{k-1} + vA^{n}(A+a_{k-1}I)zA^{k-2} + \cdots \\
& \quad  + vA^{n}(A^{k-1}+a_{k-1}A^{k-2} +\cdots + a_1I)zI,N +n].
\end{split}
\end{equation*}
We denote this homomorphism by $\varphi_{(z,N)}$.
\end{lemma}

\paragraph{Proof:}

Let $\{v_i\}_{i=1}^{K}$ be the standard basis for $\mathbb{Z}^{K}$.  For each $i$ consider $\varphi[v_i,0] = [X_i,N_i]$ where $N_i$ 
is the least integer in the equivalence class.  That is to say, if $[Y,M] = [X_i,N_i]$ then $N_i \leq M$.  Now define
$$
N = \max \{N_i \ | \ 1 \leq i \leq K \}.
$$ 
So for all $v \in \mathbb{Z}^{K}$ we can write $\varphi[v,0] = [X_v,N]$ for some $X_v \in S(A)$.  Now recall that 
$[v,n]\ast [A,1] = [v,n+1]$, and $[X_v,N] \ast [A,1] = [X_vA,N+1]$, and $\varphi$ is an $R_{A}$-module homomorphism. So
$$
\varphi[v,n] = \varphi\left( [v,0]\ast[A,1]^n \right) = \varphi[v,0] \ast [A,1]^n = [X_v,N]\otimes [A,1]^n = [X_vA^n,N+n].
$$
This shows that the number $N$ is important data in describing $\varphi$ and that $\varphi$ can be described completely by its 
restriction to terms of the form $[p(A), 0]$, where $p$ is in $\mathbb{Z}_{k-1}[x]$.
Also, for any $[v,n] \in K_0(S)$, $\varphi[v,0] $  must be of the form
$$
\varphi[v_i,0] = [(x_{k-1})_iA^{k-1} + (x_{k-2})_iA^{k-2} + \cdots + (x_{0})_iI,N],
$$
for some integers $(x_{k-1})_i, \ldots, (x_{0})_i$. It then follows from linearity that for $v \in \mathbb{Z}^{K}$ we have 
$$
\varphi[v,0] = [vx_{k-1}A^{k-1} + vx_{k-2}A^{k-2} + \cdots + vx_{0}I,N],
$$
where $x_{j} \in {\mathbb Z}^{K}$ is a column vector for each $0 \leq j \leq k-1$. 

We now use the fact that $\varphi$ is a module homomorphism to impose conditions on the $x_{i}$.  First, notice that $[v,0] \ast [A,0] = [vA,0]$, and $[X,N] \ast [A,0] = [XA,N]$.  Now, from above
$$
\varphi[vA,0] = [vAx_{k-1}A^{k-1} + vAx_{k-2}A^{k-2} + \cdots + vAx_{0}I,N]
$$
but, since $[vA,0]=[v,0]\ast[A,0]$ and $\varphi$ is a module homomorphism we have
\begin{eqnarray*}
\varphi[vA,0] &=& \varphi([v,0]\ast[A,0]) \\
&=& \varphi([v,0]) \ast [A,0] \\
&=&  [vx_{k-1}A^{k-1} + vx_{k-2}A^{k-2} + \cdots + vx_{0}I,N] \ast [A,0] \\
&=&  [vx_{k-1}A^{k} + vx_{k-2}A^{k-1} + \cdots + vx_{0}A,N] \\
&=& [vx_{k-1}(-a_{k-1}A^{k-1} - \cdots -a_0I) + vx_{k-2}A^{k-1} + \cdots + vx_{0}A,N] \\
&=& [v(x_{k-2}-a_{k-1}x_{k-1})A^{k-1} + \cdots + v(x_0-a_1x_{k-1})A -a_0vx_{k-1}I,N].
\end{eqnarray*}

Comparing coefficients of like powers of $A$ in these two expressions for $\varphi[vA,0]$, and noting that $v \in \mathbb{Z}^{K}$ was arbitrary, we see (in light of Lemma \ref{Subring-R}) that 
\begin{eqnarray*}
x_{k-2} - a_{k-1}x_{k-1} &=& Ax_{k-1} \\
x_{k-3} - a_{k-2}x_{k-1} &=& Ax_{k-2} \\
  & \vdots &  \\
x_0 - a_1x_{k-1} &=& Ax_{1} \\
-a_0x_{k-1} &=& Ax_{0}. \\
\end{eqnarray*}

So if we let $z = x_{k-1}$, then all other $x_j$'s are determined, and we have
$$
\varphi[v,0] = [vzA^{k-1} + v(A+a_{k-1}I)zA^{k-2} + \cdots + v(A^{k-1}+a_{k-1}A^{k-2} +\cdots + a_1I)zI,N].
$$

Now consider
\begin{equation*}
\begin{split}
\varphi[v,1] &= \varphi([v,0]\ast[A,1]) \\
&= [vzA^{k-1} + v(A+a_{k-1}I)zA^{k-2} + \cdots \\
& \quad + v(A^{k-1} + a_{k-1}A^{k-2} + \cdots  +  a_1I)zI,N] \ast[A,1] \\
&= [vzA^{k} \! + \! v(A+a_{k-1}I)zA^{k-1} \! + \! \cdots \! + \! v(A^{k-1}\! +\! a_{k-1}A^{k-2} \! +\! \cdots \! +\!  a_1I)zA,N+1].
\end{split}
\end{equation*}
Expanding this expression using $[A^k,N+1] = [-a_{k-1}A^{k-1} - \cdots - a_0I,N+1]$ and simplifying we are left with
\begin{equation*}
\begin{split}
\varphi[v,1] &= [vAzA^{k-1} + vA(A+a_{k-1}I)zA^{k-2} + \cdots \\
& \quad + vA(A^{k-1}+a_{k-1}A^{k-2} +\cdots + a_1I)zI,N +1],
\end{split}
\end{equation*}

Similarly, we can show inductively that
\begin{equation*}
\begin{split}
\varphi[v,n] &= [vA^{n}zA^{k-1} + vA^{n}(A+a_{k-1}I)zA^{k-2} + \cdots \\
& \quad + vA^{n}(A^{k-1}+a_{k-1}A^{k-2} +\cdots + a_1I)zI,N +n].
\end{split}
\end{equation*}
\qed

Lemma \ref{ModuleHom1} shows that each $\varphi \in Hom_{R_{A}}(K_0(S),R_{A})$ is of the form 
$\varphi_{(z,N)}$ for some $(z,N) \in \mathbb{Z}^{K} \times \mathbb{N}$.  It is also 
clear that each $(z,N) \in \mathbb{Z}^{K} \times \mathbb{N}$ gives rise to $\varphi_{(z,N)} \in Hom_{R_{A}}(K_0(S),R_{A})$ 
and that $\varphi_{(z,N)} + \varphi_{(w,N)} = \varphi_{(z+w,N)}$. It remains to be be
 determined when $\varphi_{(z,N)} = \varphi_{(w,M)}$.

\begin{lemma}\label{EquivHom}
Let $z, \ w \in \mathbb{Z}^{K}$, $N \leq M \in \mathbb{N}$. $\varphi_{(z,N)} = \varphi_{(w,M)}$ if 
and only if there exists $m \in \mathbb{N}$ such that $A^{2(m+M-N)}z = A^{2m}w$.
\end{lemma}

\paragraph{Proof:}
First suppose there exists $m \in \mathbb{N}$ such that $A^{2(m+M-N)}z = A^{2m}w$. Let $k=m+M-N$, then for all $[v,n] \in K_0(S)$
\begin{equation*}
\begin{split}
\varphi_{(z,N)}[v,n] &= [vA^{n}zA^{k-1}  + \cdots + vA^{n}(A^{k-1}+a_{k-1}A^{k-2} +\cdots + a_1I)zI,N +n] \\
&=[A^{2k}(vA^{n}zA^{k-1} + \cdots \\
& \quad + vA^{n}(A^{k-1}+a_{k-1}A^{k-2} +\cdots + a_1I)zI),M+m+n] \\
&=[A^{2m}(vA^{n}wA^{k-1} + \cdots \\
& \quad + vA^{n}(A^{k-1}+a_{k-1}A^{k-2} +\cdots + a_1I)wI),M+m+n] \\
&=[vA^{n}wA^{k-1}  + \cdots + vA^{n}(A^{k-1}+a_{k-1}A^{k-2} +\cdots + a_1I)wI ,M+n] \\
&=\varphi_{(w,M)}[v,n].
\end{split}
\end{equation*}

Now suppose $\varphi_{(z,N)} = \varphi_{(w,M)}$, for each $[v,n] \in K_0(S)$ there exists $m \in \mathbb{N}$ such that 
\begin{multline*}
A^{2(m+M-N)}(vA^{n}zA^{k-1}  + \cdots + vA^{n}(A^{k-1}+a_{k-1}A^{k-2} +\cdots + a_1I)zI) \\
= A^{2m}(vA^{n}wA^{k-1} + \cdots + vA^{n}(A^{k-1}+a_{k-1}A^{k-2} +\cdots + a_1I)wI),
\end{multline*}
However, as $l$ is the multiplicity of $0$ as a root to the minimal polynomial of $A$, $A^{l+j}X = A^{l+j}Y$ if and only if $A^lX=A^lY$, so we can replace $m$ by $l$ to get an expression which is valid for any $[v,n]$. This becomes
\begin{multline*}
[vA^{n}A^{2(l+M-N)}zA^{k-1}  + \cdots + vA^{n}(A^{k-1}+a_{k-1}A^{k-2} +\cdots + a_1I)A^{2(l+M-N)}zI,M+l+n] \\
= [vA^{n}A^{2l}wA^{k-1} + \cdots + vA^{n}(A^{k-1}+a_{k-1}A^{k-2} +\cdots + a_1I)A^{2l}wI,M+l+n].
\end{multline*}

Comparing the coefficients of $A^{k-1}$, in light of Lemma \ref{Subring-R}, we see that
$$
vA^{n}A^{2(l+M-N)}z = vA^{n}A^{2l}w
$$
for any $v \in \mathbb{Z}^{\#V(G)}$, so
$$
A^{n}A^{2(l+M-N)}z = A^{n}A^{2l}w.
$$
\qed

Thus, we have shown the following.

\begin{prop}\label{HomSys}
$Hom_{R_{A}}(K_0(S),R_{A})$ is equal to the limit of the following inductive system
$$
\xymatrix{ {\mathbb Z}^{K} \ar[r]^{z \mapsto A^2z} & {\mathbb Z}^{K} \ar[r]^{z \mapsto A^2z} & {\mathbb Z}^{K} \ar[r]^{z \mapsto A^2z} & \cdots}
$$
In other words
$$
Hom_R(K_0(S),R) \cong (\mathbb{Z}^{K} \times \mathbb{N})/ \sim_2
$$
where, for $N \leq M$,  $(z,N) \sim_2 (w,M)$ if and only if there exists $m \in \mathbb{N}$ such that $A^{2(m+M-N)}z = A^{2m}w$.
\end{prop}

Theorem \ref{ModDuality} now follows at once, since this is the same description as for $K_{0}(\USig)$.

We conclude this section with two examples.
Consider the SFT with adjacency matrix
$$
A = \left[ \begin{array}{cc}
1 & 2 \\
2 & 1
\end{array}\right].
$$

It is easy to check that $C(A)$, and hence also $K_{0}(\CH)$, are abelian, but also that $[\frac{1}{2}(A-I), 0]$ is in $K_{0}(\CH)$,
but not in $R_{A}$. Hence we have $R_{A}$ is strictly smaller than the center of $K_{0}(\CH)$.
On the other hand, with this example, we do have $Hom_{K_{0}(\CH)}(K_0(S),K_{0}(\CH))$ and  $Hom_{R_{A}}(K_0(S),R_{A})$ 
are isomorphic to each other, and to $K_0(U)$.

We now show that $Hom_Z(K_0(S),Z) \cong K_0(U)$.

Finally, consider the SFT with adjacency matrix 
$$
A = \left[\begin{array}{ccc}
0 & 1 & 5 \\
1 & 0 & 1 \\
1 & 1 & 0
\end{array}\right].
$$
Similar to the last case,  $C(A)$ and  $K_{0}(\CH)$ are abelian, but  $[\frac{1}{2}(A^{2}+A), 0]$ is in $K_{0}(\CH)$,
but not in $R_{A}$. However, it can be shown that  $Hom_{K_{0}(\CH)}(K_0(S),K_{0}(\CH))$ and  $Hom_{R_{A}}(K_0(S),R_{A})$ 
are \emph{not} isomorphic to each other.

\section{Irreducible Smale space}\label{irred}

In this section we extend the results proved in previous sections to the case that the 
Smale space is \emph{irreducible}.  The key is Smale's spectral decomposition (see \cite{ruelle78}, \cite{smale67}).

\begin{theo}[Smale]\label{irredDecompprelim}
Let $(x,\varphi)$ be an irreducible Smale space. Then there exists a positive integer $N$ and 
subsets $X_1, X_2, \ldots , X_N$ of $X$ which are closed, open, pairwise disjoint, and whose
 union equals $X$.  These sets are cyclicly permuted by $\varphi$, and $\varphi^N|_{X_i}$ is mixing for each $i$.  These sets are unique up to (cyclic) relabeling.
\end{theo}

This theorem can, alternatively, be formulated as follows.

\begin{cor}\label{irredDecomp}
Let $(X,\varphi)$ be an irreducible Smale space, then there exists a mixing Smale space $(Y,\psi)$ and a positive integer $N$ such that $X \cong Y \times \{1,\ldots,N \}$ and
$$
\varphi(y,i) = \left \{ \begin{array}{cl}
(y,i+1) & \textrm{if $1 \leq i \leq N-1$} \\
(\psi(y),1) & \textrm{if $i=N$}
\end{array} \right .
$$
\end{cor}

\paragraph{Proof:}
Let $X_1,\ldots,X_N$ be as in Prop. \ref{irredDecompprelim}.  It suffices to show that for $1\leq i \leq N-1$, $(X_i,\varphi^N) \cong (X_{i+1},\varphi^N)$ with the topological conjugacy realized by the map $\varphi$.

As $\varphi$ is a homeomorphism, it suffices to show that, for all $x \in X_i$
$$
\varphi \circ \varphi^N(x) = \varphi^N\circ \varphi(x),
$$
which is obvious. Now setting $Y=X_1$, $\psi = \varphi^N$ we have $X \cong Y \times \{1,\ldots,N \}$ and $\varphi(y,i) = (y,i+1)$ for $ 1 \leq i \leq N-1$. Finally, for all $1\leq i \leq N$, we have $\varphi^N(y,i) = (\varphi^N(y),i) = (\psi(y),i)$ which implies $\varphi(y,N) = (\psi(y),1)$.
\qed

Let $(X,\varphi)$ be a Smale space and fix $n \in \mathbb{N}$.  It is easy to see that $(X, \varphi^n)$ is also a Smale space with the same bracket function $[\cdot,\cdot]$. 
It is also easy to see that the 3 equivalence relations are unchanged by switching from $\varphi$ to $\varphi^n$.  For example, for $x \in X$ the set $V^s(x)$ is the same whether we consider the map $\varphi$ or $\varphi^n$.  In particular, for a finite $\varphi$-invariant (also $\varphi^n$-invariant) set $P \subset X$ the groupoids $\GsX$ and $G^s(X,\varphi^n,P)$ are the same.  Similarly for $\GuX$ and $G^u(X,\varphi^n,P)$, and $\GhX$ and $G^h(X,\varphi^n)$.  It then follows that $\SX = S(X,\varphi^n,P)$ and similarly for the unstable and homoclinic algebras.  It should be noted that while $\SX = S(X,\varphi^n,P)$, the automorphisms $\alpha_{\varphi}$ and $\alpha_{\varphi^n}$ are not equal.

Now suppose $(X, \varphi)$ is an irreducible Smale space and $(Y, \psi)$, $n \in \mathbb{N}$ are as in Prop. \ref{irredDecomp}.  So $(Y, \psi)$ is mixing, $X \cong Y \times \{1,2,\ldots,n \}$ and $\varphi(x,i) = (x,i+1)$ if $1\leq i \leq n-1$, $\varphi(x,n) = (\psi(x),1)$.  If we consider the Smale space $(X, \varphi^n)$ we still have $X \cong Y \times \{1,2,\ldots,n \}$, and now $\varphi^n(x,i) = (\psi(x),i)$. So $(X,\varphi^n)$ is a disjoint union of $n$ copies of the mixing Smale space $(Y,\varphi)$.

If we now fix a finite $\varphi$-invariant set $P \subset X \cong Y \times \{1,2,\ldots,n\}$, and let $\tilde{P}$ be $P \cap Y\times \{1\}$ we immediately see that
\begin{eqnarray*}
 \SX = S(X,\varphi^n,P) & \cong & \bigoplus_i^n S(Y,\psi,\tilde{P}), \\
\UX = U(X,\varphi^n,P) & \cong & \bigoplus_i^n U(Y,\psi,\tilde{P}), \ \textrm{and} \\
\HX = H(X,\varphi^n,P) & \cong & \bigoplus_i^n H(Y,\psi). \\
\end{eqnarray*}

Denote by $\alpha_{\varphi}$ and $\alpha_{\psi}$ the $\ast$-automorphisms on $\SX$ and $S(Y,\psi,\tilde{P})$ respectively. It is then straightforward to see that $\alpha_{\varphi}$ permutes the summands of $\bigoplus_i^n S(Y,\psi,\tilde{P})$. In particular, for $a \in S(Y,\psi,\tilde{P})$ we have
$$
\alpha_{\varphi}(a,i) = \left\{\begin{array}{cl}
                                (a,i+1) & 1 \leq i \leq n-1 \\
                                (\alpha_{\psi}(a),1) & i = n.
                               \end{array}\right.
$$
The corresponding results hold for $\UX$ and $\HX$ similarly.

We briefly describe the mapping cylinder and its $K$-theory ring in the case that $(X,\varphi)$ is an irreducible Smale space.

\begin{lemma}\label{IrredCH}
Let $(X,\varphi) \cong (Y,\psi) \times \{1,2,\ldots, n \}$ be an irreducible Smale space as in Prop. \ref{irredDecomp}. Then for $f \in C(H(X,\varphi),\alpha^{\varphi})$ we can write 
$$
f(t) = \left(f_1(t), f_1(t-1), \ldots, f_1(t-n+1)\right)
$$
for some $f_1 \in C(\mathbb{R},H(Y,\psi))$ such that $\alpha^{\psi}(f_1(t))=f_1(t+n)$. In other words $f_1(t) = \tilde{f}_1(\frac{t}{n})$ for some $\tilde{f}_1 \in C(H(Y,\psi),\alpha^{\psi}).$
\end{lemma}

\paragraph{Proof:}
We can write $f \in C(H(X,\varphi),\alpha^{\varphi})$ as
$$
f(t) = \left(f_1(t),f_2(t),\ldots,f_n(t)\right)
$$
where each $f_i \in C(\mathbb{R},H(Y,\psi))$ and recall that $\alpha^{\varphi}(f(t)) = f(t+1)$.  Now
$$
\alpha^{\varphi}(f(t)) = (\alpha^{\psi}(f_n(t)), f_1(t),f_2(t),\ldots,f_{n-1}(t))
$$
so, for $1\leq i \leq n-1$ we have
$$
f_i(t) = f_{i+1}(t+1), \ \textrm{or} \ f_{i+1}(t) = f_i(t-1),
$$
so for $0 \leq k \leq n-1$
$$
f_{1+k}(t) = f_1(t-k).
$$
Also,
$$
f_1(t+1) = \alpha^{\psi}(f_n(t)) = \alpha^{\psi}(f_1(t-n+1)), \ \textrm{or} \ f_1(t+n) = \alpha^{\psi}(f_1(t)).
$$
\qed

\begin{prop}\label{IrredCHIso}
 Let $(X,\varphi) \cong (Y,\psi) \times \{1,2,\ldots, n \}$ be an irreducible Smale space, then $C(H_X,\alpha_{\varphi}) \cong C(H_Y,\alpha_{\psi})$.
\end{prop}

\paragraph{Proof:}
For $f \in C(H(X,\varphi),\alpha^{\varphi})$ the map from Lemma \ref{IrredCH} which sends $f$ to $\tilde{f_1} \in  C(H(Y,\psi),\alpha^{\psi}).$ Has inverse
$$
g(t) \mapsto (g(nt), g(nt - 1), \ldots ,g(nt - n +1)).
$$
\qed

\begin{cor}\label{IrredRing}
 $K_{\ast}(C(H_X,\alpha_{\varphi})) \cong K_{\ast}(C(H_Y,\alpha_{\psi}))$ as rings.
\end{cor}

\begin{remark}\label{IrredModRmk}
It is now straightforward to write down the module structures for an irreducible Smale space. In particular, if $(X,\varphi)$ is an irreducible Smale space then
$$
\SX \cong \bigoplus_1^n S(Y,\psi,\tilde{P})
$$
and
$$
C(H_X,\alpha_{\varphi}) \cong C(H_Y, \alpha_{\psi}).
$$
So for $[b] \in K_{\ast}C(H_Y, \alpha_{\psi})$, $([a_1],[a_2],\ldots,[a_n]) \in \SX$ we have
$$
([a_1],[a_2],\ldots,[a_n])[b] = ([a_1][b],[a_2][b],\ldots,[a_n][b]).
$$
The corresponding results for the module structures on $\UX$ and $\HX$ also hold.

\end{remark}

\end{document}